%% file: SLLN_immigration_26.tex
\documentclass[11pt]{article}
\usepackage{latexsym,amsmath,amssymb,amsfonts,amsbsy,fancyhdr}
\usepackage{theorem}
\usepackage{amscd}
\usepackage{enumerate}
\usepackage{mathrsfs}
\usepackage{dsfont}
\usepackage[active]{srcltx}
\usepackage[colorlinks=true, pdfstartview=FitV, linkcolor=black, citecolor=black, urlcolor=black]{hyperref}
\usepackage[dvips]{graphicx}
\usepackage{epsf}
\usepackage{nicefrac}
\usepackage{marvosym}
\usepackage{subfigure}

\newcommand{\R}{{\mathbb R}}

\newcommand{\rpn}{\R^+_0}

\newcommand{\N}{{\mathbb N}}

\newcommand{\iid}{\text{i.i.d.} }

\newcommand{\dd}{\text{d}}

\newtheorem{theorem}{Theorem}
\newtheorem{proposition}{Proposition}
\newtheorem{corollary}[proposition]{Corollary}
\newtheorem{lemma}[proposition]{Lemma}

\theorembodyfont{\rm}
\newtheorem{definition}[proposition]{Definition}
\newtheorem{example}[proposition]{Example}

\newtheorem{remark}[proposition]{Remark}
\newenvironment{rem}{\begin{remark}}{\hfill$\lozenge$\end{remark}}
\newtheorem{proof}{Proof}

\newenvironment{pf}{\begin{proof}}{\hfill$\square$\end{proof}}
\newtheorem{hypothesis}{Hypothesis}

\allowdisplaybreaks

\DeclareMathAlphabet{\mathcal}{OMS}{cmsy}{m}{n}

\begin{document}

\title{\Large \bf Limit theorems for fragmentation processes with immigration}
\author{Robert Knobloch\thanks{Institut f\"ur Mathematik, FB 12, Goethe-Universit\"at Frankfurt am Main, 60054 Frankfurt am Main, Germany
\newline
e-mail: knobloch@math.uni-frankfurt.de}}
\date{\today}
\maketitle

\begin{abstract}
In this paper we extend two limit theorems which were recently obtained for fragmentation processes to such processes with immigration. More precisely, in the setting with immigration  we consider the asymptotic behaviour of an empirical measure associated with the stopping line corresponding to the first blocks, in their respective line of descent, of size less than $\eta\in(0,1]$ as well as a limit theorem for the process counted with a random characteristic. In addition, we determine the asymptotic decay rate of the size of the largest block in a homogeneous fragmentation process with immigration. The techniques used to proves these results are based on submartingale arguments.
\end{abstract}



\section{Introduction}
The present paper deals with some asymptotic properties  of fragmentation processes with immigration. Random fragmentations with immigrations were introduced in \cite{98}, where also the deterministic setting was considered. The model this paper is concerned with consists of the interplay of two independent random actions with opposing effects. On the one hand there is the fragmentation of the blocks present in the system, which increases the number of blocks but may decrease the total mass in the system. On the other hand, new blocks immigrate into the system and  increase the number of blocks as well as the total size of all the blocks present in the system. In \cite{98} Haas investigated the existence and uniqueness of a stationary distribution for fragmentation  processes with immigration. 
In this model the immigration and the fragmentation events are described  by independent  Poisson point processes. More precisely, at each time at which the immigration Poisson point process has an atom at most countably many  fragmentation processes immigrate into the system. Each immigrated block then evolves as a fragmentation process (cf. Section~\ref{s.hfp}), that is it  fragments into a collection of smaller block according to a Poisson point processes which is independent of the immigration Poisson point processes and independent of the fragmentation of the other blocks.  Note that the total number of blocks remains at most countably infinite but the total size may be unbounded. In particular, by not taking the times of immigration into account, fragmentations with immigration can be considered as fragmentation processes whose total mass is finite at every time but not bounded by a  constant. 

In this paper we extend the main results of \cite{HKK10} and \cite{Kno12}  to fragmentation processes with immigration. That is to say, the goal we are concerned with is tantamount to showing that in a certain manner the limiting behaviour of the whole system is the same as taking the limit of the object under consideration for each immigrating block separately and putting these limits together. 

In the setting with immigration we prove for bounded, deterministic functions the  almost sure convergence of the empirical mean that is defined via a fragmentation process with immigration which is stopped at the stopping line of first passage times below a given block-size. In the setting of fragmentation chains (without immigration) this kind of problem was considered in \cite{BM05} as well as \cite{HK11} and we shall resort to \cite{HKK10} where the corresponding strong law of large numbers was obtained for general fragmentation processes without immigration.

Moreover,  we prove the convergence of the so-called process counted with a random characteristic. Random characteristics were considered by various authors for general branching processes and in the context of fragmentation processes they were introduced in \cite{Kno12}. As for applications of such processes counted with a random characteristic for fragmentation processes we refer for instance to \cite{BM05}, where convergence results for the fragmentation energy model  were established.

A similar motivation as for the abovementioned extensions to the setting with immigration lies at the heart of Theorem~4.2 in \cite{Olo96}, where Olofsson extends an $\mathscr L^1$-convergence result that was obtained for Crump-Mode-Jagers processes in \cite{Ner81} to such processes with immigration. We also refer to \cite{99} for a result on supercritical immigration-branching processes  in the spirit of this paper. 


The final result of this paper is concerned with the asymptotic decay of the largest block in the process with immigration. We give bounds on the exponential decay rate, where the upper bound depends on the time at which large blocks immigrate.

The outline of this paper is as follows. In the next section we introduce  fragmentation processes and such processes stopped at a particular family of stopping lines. The third section aims at providing some  motivation for the topic of this paper by describing the problems under consideration in the setting without immigration. Subsequently, in Section~\ref{s.setup} we establish the set-up for fragmentations with immigration and state the main results of this paper. Section~\ref{s.example} deals with an example of a fragmentation process with immigration that is based on the spine-decomposition of a fragmentation process. The last three sections are concerned with the proofs of our main results.

Throughout this paper we consider a probability space $(\Omega,\mathscr F,\mathbb P)$ on which the fragmentation processes as well as all the other random objects are defined.

\section{Stopped fragmentation processes}\label{s.hfp}
In this section we provide a brief introduction to fragmentation processes and associated  processes that are stopped at certain stopping lines. Moreover, we introduce two families of additive martingales that we shall need later on. 

We shall consider two closely related classes of fragmentation processes, namely partition-valued fragmentations and mass fragmentations, respectively. Let us start by defining the state space of the partition-valued processes. To this end, let $\mathcal{P}$ be the space of partitions $\pi=(\pi_n)_{n\in\N}$ of $\N$, where the blocks of $\pi$ are ordered by their least element such that $\inf(\pi_i)<\inf(\pi_j)$ if $i<j$, where $\inf(\emptyset):=\infty$. This paper is concerned with a $\mathcal P$-valued fragmentation process $\Pi:= (\Pi(t))_{\eta\in(0,1]}$, where $\Pi(t) = (\Pi_n(t))_{n\in\N}$. $\mathcal{P}$-valued fragmentations are exchangeable Markov processes which were introduced in  \cite{85} in the homogeneous case and were extended to the self-similar setting in \cite{Ber02b}. For a comprehensive treatise on fragmentation processes we refer to the monograph  \cite{Ber06}. Let $\mathscr F:=(\mathscr F_t)_{\eta\in(0,1]}$ be the filtration generated by $\Pi$. 

It is known from \cite{Ber02b} that the distribution of $\Pi$ is determined by some $\alpha\in\R$ (the index of self-similarity; $\alpha=0$ corresponding to the homogeneous case), a constant $c\in\rpn$ (the rate of erosion) and a measure $\nu$ (the so-called dislocation measure that determines the jumps of $\Pi$) on 
 \[
\mathcal S_1:=\left\{{\bf s}:=(s_n)_{n\in\N}:s_1\ge s_2\ldots\ge0,\,\sum_{n\in\N}s_n\le1\right\},
 \] 
such that  $\nu(\{(1,0,\cdots)\})=0$ as well as 
\[
\int_{\mathcal S_1}(1-s_1)\nu(d{\bf s})<\infty.
\]
The measure $\nu$ is said to be {\it conservative} if $\nu(\sum_{n\in\N}s_n<1)=0$, i.e. if there is no loss of mass by sudden dislocations, and {\it dissipative} otherwise. In this paper we allow for both of these cases. 
Below we shall need the following constant 
\[
\underline{ p}: = \inf\left\{  p\in \mathbb{R} : \int_{\mathcal S_1} \left| 1- \sum_{n\in\N} s_n^ {1+p} \right| \nu(d{\bf s}) <\infty\right\}\in(-1, 0]
\]
as well as the increasing and concave function $\Phi:(\underline p,\infty)\to\R$, given by
\[
\Phi(p)=\int_{\mathcal S_1}\left(1-\sum_{n\in\N}s_n^{1+p}\right)\nu(\dd{\bf s})
\]
for every $(\underline p,\infty)$. If $\underline p=0$ in the conservative case, then we set $\Phi(\underline p):=0$.  The function $\Phi$ plays a crucial role in the theory of fragmentation processes, since it turns out to be the Laplace exponent of the killed subordinator $\xi:=(\xi_t)_{t\in\rpn}$ defined by
\[
\xi_t:=-\ln(|\Pi_1(t)|)
\]
for all $t\in\rpn$. 

{\it Throughout this paper, unless stated otherwise, we consider a homogeneous fragmentation process $\Pi$  that satisfies $c=0$ as well as the two hypotheses which we are just about to introduce.} On this note let us point out that, as outlined in Corollary~\ref{c.self.sim}, many of the results of this paper also hold for self-similar fragmentation processes.

If there exists a $p^*\in[\underline p,0]$\label{p.p^*} satisfying $\Phi(p^*)=0$, then we call $p^*$ {\it Malthusian parameter}. The following Hypothesis~\ref{h.m}, commonly referred to as {\it Malthusian hypothesis}\index{Malthusian hypothesis}, provides us with the existence of a Malthusian parameter in the dissipative case. 

\begin{hypothesis}\label{h.m}
If $\Pi$ is dissipative, then there exists a $p^*\in(\underline p,0)$ such that $\Phi(p^*)=0$. 
\end{hypothesis}

If $\Pi$ is conservative, that is if $\nu\left(\sum_{n\in\N}s_n<1\right)=0$, then $\Phi(0)=0$ and thus we set $p^*:=0$ in that case. In view of Lemma~1 in \cite{Ber03} let $\bar p$ be the unique solution to 
\[
(1+p)\Phi'(p)=\Phi(p)
\]
on $(\underline p,\infty)$, where $\Phi'$\label{p.Phi'} denotes the derivative of $\Phi$.   Notice that it follows from Lemma~1 in \cite{Ber03} that $p\ge\bar p$ if and only if $(1+p)\Phi'(p)\le\Phi(p)$. Since $\Phi'(p)>0$ for all $p\in(\underline p,\infty)$, we therefore have $p^*<\bar p$. Moreover, observe that Hypothesis~\ref{h.m} implies that $\underline p<0$ and thus $\Phi'(0+)<\infty$ in the dissipative case. However, in the conservative case it is possible that $\underline p=0$, in which case the expectation of the subordinator $\xi$ may be infinite. In order to guarantee that $\xi$ has finite expectation in the conservative case, we need the following hypothesis:

\begin{hypothesis}\label{h.a}
If $\underline p=0$, then 
\[
\Phi'(0+)=\int_{\mathcal S_1}\left(\sum_{n\in\N}s_n\ln\left(s_n^{-1}\right)\right)\nu(\dd\text{s})<\infty.
\]
\end{hypothesis}
Henceforth, we assume that Hypothesis~\ref{h.m} and Hypothesis~\ref{h.a} hold. 


We shall need the  exchangeable partition measure $\mu$ on $\mathcal{P}$ given by
\begin{equation}\label{e.mu}
\mu(d\pi) = \int_{\mathcal S_1}\varrho_{\bf s}(d\pi)\nu(d{\bf s}),
\end{equation}
where $\varrho_{\bf s}$ is the law of Kingman's paint-box based on ${\bf s}\in\mathcal S_1$. In \cite{85} Bertoin showed that the homogeneous  fragmentation process $\Pi$ is characterised by a Poisson point process. More precisely, there exists a $\mathcal P\times\mathbb N$-valued Poisson point process $(\pi(t),\kappa(t))_{t\in\rpn}$\label{p.pi_t} with characteristic measure $\mu\otimes\sharp$, where $\sharp$ denotes the counting measure on $\N$,  such that $\Pi$ changes state precisely at the times $t\in\rpn$ for which an atom $(\pi(t),\kappa(t))$ occurs in $(\mathcal P\setminus(\N,\emptyset,\ldots))\times\mathbb N$. At such a time $t\in\rpn$ the sequence $\Pi(t)$ is obtained from $\Pi(t-)$ by replacing its $\kappa(t)$-th term, $\Pi_{\kappa(t)}(t-)\subseteq\N$, with the restricted partition $\pi(t)|_{\Pi_{\kappa(t)}(t-)}$ and reordering the terms such that the resulting partition of $\N$ is an element of $\mathcal P$. We denote the random jump times of $\Pi$, i.e. the times at which the abovementioned Poisson point process has an atom in $(\mathcal P\setminus(\N,\emptyset,\ldots))\times\mathbb N$, by $(t_i)_{i\in\mathcal I}$, where the index set $\mathcal I\subseteq\rpn$ is  countably infinite.  

Moreover, by exchangeability the limits
\[
|\Pi_n(t)|:=\lim_{k\to\infty}\frac{\sharp(\Pi_n(t)\cap\{1,\ldots,k\})}{k},
\]
referred to as {\it asymptotic frequencies}, exist $\mathbb P$-a.s. simultaneously for all $t\in\rpn$ and $n\in\N$. Let us point out that the concept of asymptotic frequencies provides us with a notion of {\it size} for the blocks of a $\mathcal P$-valued fragmentation process. We denote by $\lambda:=(\lambda(t))_{t\in\rpn}$, with $\lambda(t)=(\lambda_n(t))_{n\in\N}$, the $\mathcal S_1$-valued {\it mass fragmentation process} given by
\begin{equation}\label{e.lambda}
\lambda_n(t)=|\Pi(t)|^\downarrow_n
\end{equation}
for all $n\in\N$ and $t\in\rpn$, where $(|\Pi(t)|^\downarrow_n)_{n\in\N}$ denotes the decreasing reordering of the asymptotic frequencies of the blocks of $\Pi(t)$. In addition, let $\mathscr G:=(\mathscr G_t)_{t\in\rpn}$ be the filtration generated by $\lambda$. In this respect we also refer to \cite{Ber02a}, where mass fragmentation processes (also known as ranked fragmentations) were considered.

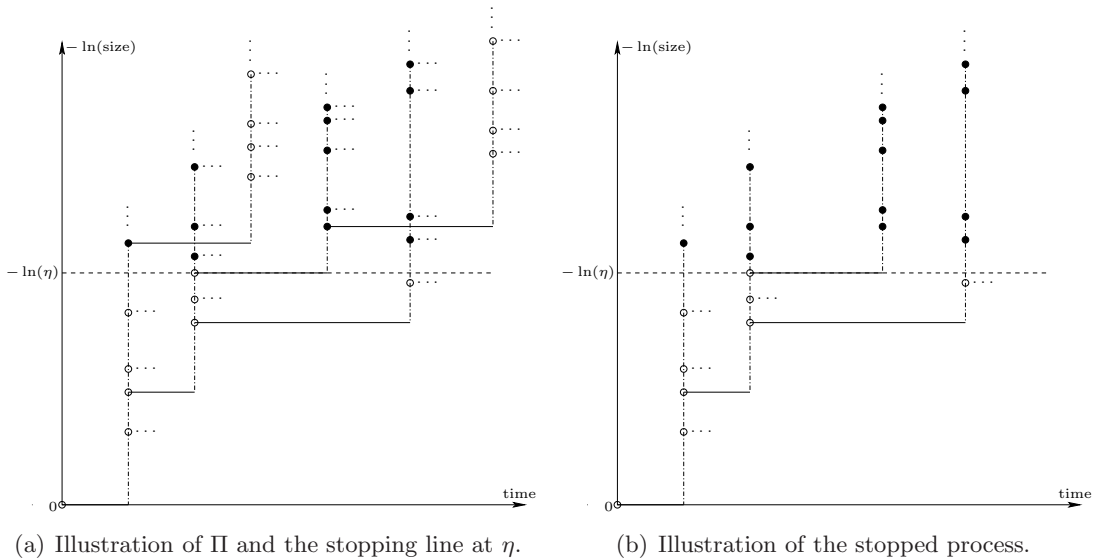
\begin{figure}[htb]\centering
\subfigure[Illustration of $\Pi$ and the stopping line at $\eta$.\label{f.1.2a}]{\resizebox{7cm}{!}{\input{frag_stopping_line_50_eta.pstex_t}}}\quad
\subfigure[Illustration of the  stopped process.\label{f.1.2b}]{\resizebox{7cm}{!}{\input{frag_stopping_line_50b_eta.pstex_t}}}\quad
\caption[Stopped fragmentation process]{Illustration (a) depicts a realisation of a  fragmentation process $\Pi$ (with finite dislocation measure) and the stopping line  given by the first passage of the block sizes below some $\eta\in(0,1]$. In (b)  the fragmentation process  which is stopped at this stopping line is illustrated. The black dots indicate the blocks in $\lambda_\eta$, since their sizes are smaller than $\eta$ and they result from the dislocation of blocks with size greater than or equal to $\eta$. 
}\label{f.1.2}
\end{figure}

The main  property of fragmentation processes 
 is the {\it (strong) fragmentation property}, which is the analogue of the branching property of branching processes. Roughly speaking, this property says that given a configuration of the process at some (stopping) time, the  further evolution of each block is governed by an independent copy of the original process. Moreover, the same holds true if we replace the stopping time by a stopping line, cf. Definition~3.4 and Lemma~3.14 both in \cite{Ber06}, in which case we refer to it as the {\it extended fragmentation property}. Speaking of stopping lines, let us now introduce stopped fragmentations which are obtained from a fragmentation process by stopping the evolution of a block once it has reached a given stopping line. Here we are interested in the stopping line that corresponds to the first passage times
\begin{equation}\label{e.stopping_line}
\upsilon_{\eta,k}:=\inf\left\{t\in\rpn:|B_k(t)|<\eta\right\},\quad k\in\N,\,\eta\in(0,1], 
\end{equation} 
i.e. to the first blocks, in their respective line of descent, of size less than some given $\eta\in(0,1]$. The  process $(\lambda_\eta)_{\eta\in(0,1]}$ consisting for each $\eta\in(0,1]$ of the sizes of the blocks at the terminal state of the fragmentation  stopped at the stopping line associated with $\eta$   is then given by $\lambda_\eta:=(\lambda_{\eta,k})_{k\in\N}$, where $\lambda_{\eta,k}$ refers  to the asymptotic frequency of the $k$-th largest block at the terminal state of the stopped process. In addition, we denote by $(\mathscr H_\eta)_{\eta\in(0,1]}$ the filtration generated by $(\lambda_\eta)_{\eta\in(0,1]}$. For an illustration of 
these concepts, see Figure~\ref{f.1.2}.
\\

In what follows we shall make use of two 
unit-mean martingales. The first martingale, $M(p):=(M_t(p))_{t\in\rpn}$, is given by
\begin{equation}\label{e.martingale_M}
M_t(p):=\sum_{k\in\N}|\Pi_{k}(t)|^{1+p}e^{\Phi(p)t}
\end{equation}
for every $p\in(\underline p,\bar p)$ and $t\in\rpn$. Let us point out that $M(p)$ is the analogue of Biggins'  additive martingale for branching random walks (cf. \cite{Big77}). These martingales  appear frequently in the literature on branching processes and fragmentation processes. In addition,  we shall need a second  martingale which is defined in terms of $(\lambda_\eta)_{\eta\in(0,1]}$ and turns out to be related to $M(p^*)$. More precisely, consider the process $\Lambda(p^*):=(\Lambda_\eta(p^*))_{\eta\in(0,1]}$   defined by
\[
\Lambda_\eta(p^*):=\sum_{k\in\N}\lambda_{\eta,k}^{1+p}
\]
for each $\eta\in(0,1]$
The martingale property of $\Lambda(p^*)$ was established in Lemma~1 of \cite{HKK10}. In particular, there it was shown that 
\[
\Lambda_\eta(p^*)=\mathbb E\left(\left.M_\infty(p^*)\right|\mathscr H_\eta\right)
\]
holds for all $\eta\in(0,1]$, where $M_\infty(p^*):=\lim_{t\to\infty}M_t(p^*)$. As we will see in the next section, the limits which  appear in our main results turn out to be deterministic constants times the almost sure martingale limit 
\begin{equation}\label{e.Lambda_0}
\Lambda_0(p^*):=\lim_{\eta\downarrow0}\Lambda_\eta(p^*).
\end{equation}

\section{Motivation}\label{s.motivation}
The principal goal of this paper is to extend two results which were recently obtained for  fragmentation processes to the setting of fragmentations with immigration. In this spirit  let us devote this section to a description in the setting without immigration of the results that we aim at extending.

\subsection{Empirical distribution for stopped fragmentations}\label{s.edsf}
The first problem which we shall be concerned with in the context of fragmentation processes with immigration deals with the asymptotic behaviour of an empirical measure associated with a certain family of stopping lines. Before we state the corresponding result in the setting without immigration, let us first establish some notation. For this purpose, let $\mathcal B^+_b$\label{p.B+} denote the space of all bounded and measurable functions $f:\rpn\to\rpn$ with $f|_{[1,\infty)}\equiv0$. 
Moreover, for any $\eta\in(0,1]$ consider the random measure $\rho_\eta$\label{p.rho_t} on $[0,1]$ defined by
\[
\rho_\eta:=\sum_{k\in\N} \lambda_{\eta,k}^{1+p^*} \delta_{\frac{\lambda_{\eta,k}}{\eta}},
\]
where $\delta_x$\label{p.delta_x} is the Dirac measure at $x\in[0,1]$. 
The result we  present in this section  is concerned with the associated empirical mean, i.e. with the integral of test functions in $\mathcal B^+_b$ with respect to the empirical measure $\rho_\eta$. In this regard we define
\begin{equation}\label{e.em.1}
\langle\rho_\eta,f\rangle:=\int_{[0,1]}f\text{ d}\rho_\eta=\sum_{k\in\N}\lambda_{\eta,k}^{1+p^*} f\left(\frac{\lambda_{\eta,k}}{\eta}\right)
\end{equation}
for any $\eta\in(0,1]$ and $f\in\mathcal B^+_b$. Notice that the time-parameter of the process $(\langle\rho_\eta,f\rangle)_{\eta\in(0,1]}$ corresponds to the size rather than to the time of the fragmentation process $\Pi$. Define a measure $\rho$\label{p.rho} on $[0,1]$ as follows:
\[
\rho(\dd t):=\frac{1}{\Phi'(p^*)}\left(\int_{\mathcal S_1}\sum_{n\in\N}\mathds1_{\{s_n<t\}}s_n^{1+p^*}\nu(\dd{\bf s})\right)\frac{\dd t}{t},
\]
where in the case $p^*=\underline p=0$ we adopt $\Phi'(p^*)=\Phi'(0+)$. 
In Theorem~1 of \cite{HKK10} it was shown  that asymptotically as $\eta\downarrow0$ the random function $\eta\mapsto\langle\rho_\eta,\cdot\rangle$ behaves $\mathbb P$-a.s. like the limit $\Lambda_0(p^*)$ of the nonnegative martingale $\Lambda(p^*)$, up to a multiplicative function $\langle\rho,\cdot\rangle$ on $\mathcal B^+_b$ given by
\[
\langle\rho,f\rangle:=\int_{(0,1)}f\dd\rho=\frac{1}{\Phi'(p^*)}\int_{(0,1)}f(t)\left(\int_{\mathcal S_1}\sum_{n\in\N}\mathds1_{\{s_n<t\}}s_n^{1+p^*}\nu(\dd{\bf s})\right)\frac{\dd t}{t}
\]
for all $f\in\mathcal B^+_b$. Resorting to Theorem~1 of \cite{BHS99} and applying Lemma~2 of \cite{HKK10} one can show that 
\begin{equation}\label{p.rho_f_b}
\langle\rho,f\rangle=\lim_{\eta\downarrow0}\mathbb E(\langle\rho_\eta,f\rangle).
\end{equation}
The goal of this paper is to extend the abovmentioned  result that was proven in \cite{HKK10}. In this respect see also Theorem~4 of \cite{Kno12}.
\begin{proposition}[Theorem~1 of \cite{HKK10}]\label{t.1.1}
For any $f\in\mathcal B^+_b$ we have
\begin{equation}\label{e.convergence_em}
\langle\rho_\eta,f\rangle\to\Lambda_0(p^*)\langle\rho,f\rangle
\end{equation}
$\mathbb P$-a.s. as $\eta\downarrow0$.
\end{proposition}
Recall that we assume $\Pi$ to satisfy Hypothesis~\ref{h.m} as well as Hypothesis~\ref{h.a}. On this note let us mention that Theorem~1 of \cite{HKK10} was proven under an additional third assumption (cf. (A3) in \cite{HKK10}). However, it follows by means of Proposition~3.5 in \cite{Kno11} (cf. also Theorem~2 in \cite{Ber03} for the conservative case) that also without this additional assumption the proof provided in \cite{HKK10} works.

Proposition~\ref{t.1.1} follows a series of strong laws of large numbers that were obtained for different classes of stochastic processes with a branching structure. Related classical strong laws were considered in \cite{AH76} and \cite{99} for spatial branching processes. Nerman \cite{Ner81} proved a more general strong law of large numbers in the setting of general branching processes, from which the conclusion of Proposition~\ref{t.1.1} follows for that class of branching processes. We consider such a general result in Section~\ref{ss.rc}, but in the context of fragmentation processes the connection with Proposition~\ref{t.1.1} is not as clear as for general branching processes, see Theorem~4 of \cite{Kno12}  in this regard. In addition, strong laws of large numbers in the spirit of Proposition~\ref{t.1.1} were obtained in \cite{CS07}, \cite{CRW08} and \cite{EHK10} for branching diffusions and in \cite{EW06} as well as \cite{Eng09} in the setting of superdiffusions. For related results, in the sense of $\mathscr L^2$-convergence, in the context of conservative fragmentation processes we refer to  \cite{BM05} as well as \cite{HK11}, where the latter is concerned with conservative fragmentation chains. 
An extension of Proposition~\ref{t.1.1} to certain random functions $f$  was given in Theorem~4 of \cite{Kno12} and thus our forthcoming Theorem~\ref{t.1.2} also extends to such random $f$ by resorting in its proof to Theorem~4 of \cite{Kno12} instead of applying Proposition~\ref{t.1.1}.

\subsection{Random characteristics for fragmentations}\label{ss.rc}
The result that we describe here is in spirit similar to the one described in the previous section. Indeed, in a sense it is more general and can in fact be used to prove Proposition~\ref{t.1.1}, see the proof of Theorem~4 in \cite{Kno12}.

In \cite{Kno12} a  limit theorem for the process counted with a random characteristic was proven in the context of self-similar fragmentation processes. Here we aim at extending this result. To this end  we define a {\it random characteristic} as a  random function $\phi:\rpn\times\mathcal P\times\Omega\to\R$ such that $(\phi(x,\pi))_{x\in\rpn}$, $\pi\in\mathcal P$, is an $\R$-valued  stochastic processes,  with $\phi(x,\pi)=\phi(0,\pi)=0$ $\mathbb P$-a.s. for all $x>1$, which 
has c\`adl\`ag paths $\mathbb P$-a.s. and is independent of $\Pi$. In addition, we assume that $\phi(\cdot,(1,0,\ldots))=0$.

The object under consideration in \cite{Kno12} is the process $(Z^\phi_\eta)_{\eta\in(0,1]}$, often referred to as the {\it process counted with the characteristic $\phi$}, 
given by
\[
Z^\phi_\eta=\sum_{i\in\mathcal I}\phi_i\left(\frac{\eta}{|\Pi_{\kappa(t_i)}(t_i-)|},\pi(t_i)\right)\mathds1_{\left\{|\Pi_{\kappa(t_i)}(t_i-)|\ne0\right\}}
\]
for each $\eta\in(0,1]$, where the $\phi_i$ are independent copies of $\phi$. Interesting examples of such processes counted with a random characteristics in the context of fragmentation processes are for instance the fragmentation energy needed to crush blocks until all the fragments are smaller than a given threshold (see e.g. \cite{BM05}) as well as the empirical mean with which Section~\ref{s.edsf} is concerned.

The result that we aim at extending to the setting with immigration is concerned with the asymptotic behaviour of $Z^\phi_\eta$ as $\eta\downarrow0$. 
\begin{proposition}[Theorem~1 of \cite{Kno12}]\label{t.rc.1}
Assume that $\Pi$  satisfies
\begin{equation}\label{e.assumption_t.rc}
\nu(\mathcal S_1)=\infty\qquad\text{and}\qquad\nu({\bf s}\in\mathcal S_1:s_2=0)<\infty\qquad\text{as well as}\qquad p^*>\underline p.
\end{equation}
Further, let $\phi$ be a random characteristic such that 
\begin{equation}\label{e.assumption_1}
\int_{\mathcal P}\mathbb E\left(\sup_{\eta\in(0,1]}\eta^{(1+p^*+\beta)}\phi\left(\eta,\pi\right)\right)\mu(\dd\pi)<\infty\quad\text{and}\quad\limsup_{\eta\downarrow0}\int_{\mathcal P}\eta^{1+\tilde p}\mathbb E\left(\phi\left(\eta,\pi\right)\right)\mu(\dd\pi)<\infty
\end{equation}
hold for  all $\beta>0$ and some $\tilde p\in(\underline p,p^*)$.
Then 
\begin{equation}\label{e.convergence_rc}
\eta^{1+p^*}Z^\phi_\eta\to \frac{\Lambda_0(p^*)}{\Phi'(p^*)}\int_{(0,1)}\mathbb E\left(\sum_{k\in\N}|\Pi_{k}(t)|^{1+p^*}\int_\mathcal P\int_{(0,1]}u^{p^*}\phi\left(u,\pi\right)\dd u\mu(\dd\pi)\right)\dd t
\end{equation}
$\mathbb P$-a.s. and in $\mathscr L^1(\mathbb P)$ as $\eta\downarrow0$.
\end{proposition}

Similar results for general branching processes were obtained in \cite{Ner81} and the $\mathscr L^1$-convergence was extended to such processes with immigration in \cite{Olo96}. As mentioned above, an example of a process counted with a random characteristic, which is meaningful also from an application-oriented point of view, describes the total energy cost needed to crush blocks of unit size into blocks of size less than some given $\eta\in(0,1]$. Convergence results for the process counted with this particular characteristic were obtained in \cite{BM05} as well as \cite{Kno12}. In this regard we also refer to \cite{FKM10}. 

Notice that up to the deterministic constants the random object in the limits of (\ref{e.convergence_em}) and (\ref{e.convergence_rc}) coincide. In fact, this is no coincidence, since in principle the convergence in (\ref{e.convergence_em}) could be considered as a special case of (\ref{e.convergence_rc}), if one could show that the assumptions of Proposition~\ref{t.rc.1} are satisfied. However, even though the statement of Proposition~\ref{t.1.1} can be proven by means of Proposition~\ref{t.rc.1}, it is not as straightforward as for the analogous results in the setting of  general branching processes.


\section{Set-up and main results}\label{s.setup}
In this section we introduce self-similar fragmentation processes with immigration and we present the main results of this paper. 

\subsection{Fragmentation processes with immigration}\label{intro_immigration}
Here we introduce the setting of immigration for fragmentation processes. In this respect, see also \cite{98} for a similar set-up. Let us start by considering the state space of the process with immigration. For this purpose, set
\[
\mathcal S:=\left\{{\bf s}:=(s_n)_{n\in\N}:\sum_{n\in\N}s_n<\infty,\,0\le s_j\le s_i<\infty\,\forall\,i\le j\right\}.\label{p.S}
\]

On $\mathcal S$ we define the binary operator $+$ as the decreasingly ordered concatenation of two sequences in $\mathcal S$. The corresponding iterated operator is denoted by $\sum$.

Let ${\bf u}:=(u_n)_{n\in\N}:\Omega\to\mathcal S$. Then we call {\it self-similar mass fragmentation process starting from ${\bf u}$}\index{fragmentation process!starting from ${\bf u}$} with index $\alpha_{\bf u}:=(\alpha_n)_{n\in\N}$, $\alpha_n\in\R$, the $\mathcal S$-valued Markov process $\lambda^{\bf u}:=(\lambda^{\bf u}(t))_{\eta\in(0,1]}$\label{p.lambda^u}, defined by
\[
\lambda^{\bf u}(t):=\sum_{n\in\N}u_n\lambda^{(n)}(u^{\alpha_n}_nt)
\]
for all $t\in\rpn$, where the $\lambda^{(n)}$ are independent 
self-similar  mass fragmentation processes with index $\alpha_n$, where we assume that the $\lambda^{(n)}$ are also independent of ${\bf u}$.

Let us now define fragmentation processes with immigration.

\begin{definition}\label{d.fpi.2}
Consider a Poisson measure $N$ on $\mathcal S\otimes\rpn$ with intensity $\nu_I\otimes\dd t$, where $\nu_I$ is a  $\sigma$-finite measure on $\mathcal S$, satisfying
\begin{equation}\label{e.immigration_measure}
\int_\mathcal S\sum_{k\in\N}(s_k\land1)\nu_I({\bf s})<\infty,
\end{equation}
and $\dd t$ is the Lebesgue measure on $\rpn$. Denote the  atoms in $\mathcal S\setminus\{(0,\ldots)\}$ of the associated Poisson point process by $({\bf s}(r_n))_{n\in\mathcal N}$, where the index set $\mathcal N$ is at most countably infinite, 
and let ${\bf u}:\Omega\to\mathcal S$ be independent of $N$. Furthermore, conditionally on {\bf u} and $({\bf s}(r_n))_{n\in\mathcal N}$ let $\lambda^{\bf u}$  and $\lambda^{{\bf s}(r_n)}$, $n\in\mathcal N$, be independent 
self-similar mass fragmentation processes starting from ${\bf u}$ and ${{\bf s}(r_n)}$ respectively. 
Then we call the $\mathcal S$-valued process $\lambda^I:=(\lambda^I(t))_{t\in\rpn}$\label{p.lambda^I}, defined by
\[
\lambda^I(t):=\lambda^{\bf u}(t)+\sum_{n\in\mathcal N:r_n\le t}\lambda^{{\bf s}(r_n)}(t-r_n)
\]
for all $t\in\rpn$, a {\it self-similar mass fragmentation process with immigration}\index{fragmentation process!\midtilde\,\,with immigration} starting from ${\bf u}$.
\end{definition}


Notice that (\ref{e.immigration_measure}) implies that $\lambda^I(t)\in\mathcal S$ $\mathbb P$--a.s. for every $t\in\rpn$. Note further  that in Definition~\ref{d.fpi.2}  the intensity of the Poisson measure $N$ is of the form $\nu_I\otimes\dd t$, i.e. the distribution of the immigrating particles can depend on the time at which the particles are immigrating.

Let ${\bf u}$ and $({\bf s}(r_n))_{n\in\mathcal N}$ be as in Definition~\ref{d.fpi.2}.  In view of the countability of $\mathcal N$ let $\mathcal J$ be an at most  countably infinite index set  and let $(v_j)_{j\in\mathcal J}$ be a vector which consists of ${\bf u}$ and $({\bf s}(r_n))_{n\in\mathcal N}$. That is to say, for every $j\in\mathcal J$ there exists some $n\in\N$ or $(n,k)\in\mathcal N\times\N$ such that $v_j=u_n$ or $v_j={\bf s}_k(r_n)$ and vice versa for every $n\in\N$ and $(m,k)\in\mathcal N\times\N$ there exist $j_n,j_{m,k}\in\mathcal J$ such that $v_{j_n}=u_n$ and $v_{j_{m,k}}={\bf s}_k(r_m)$.  
For every  $j\in\mathcal J$ let $\Pi^{(j)}$  be a homogeneous  $\mathcal P$-fragmentation  process, which we associate with $v_j$, and let $\lambda^{(j)}$ be the associated mass fragmentation process as defined in (\ref{e.lambda}). 
Furthermore, for every $j\in\mathcal J$ we denote the time at which $\Pi^{(j)}$ immigrates by $\mathfrak t_j$. Notice for any $j\in\mathcal J$ that $\mathfrak t_j=0$ if $v_j=u_n$ for some $n\in\N$ and that $\mathfrak t_j=r_n$ if $v_j={\bf s}_k(r_n)$ for some $n\in\mathcal N$ and $k\in\N$.  We assume  that the $\Pi^{(j)}$, $j\in\mathcal J$,  are independent and that $\Pi^{(j)}$ is independent of  $(v_j)_{j\in\mathcal J}$ and $(\mathfrak t_j)_{j\in\mathcal J}$. Moreover, in the light of Section~\ref{s.hfp} the $\Pi^{(j)}$, $j\in\mathcal J$, are assumed to satisfy Hypothesis~\ref{h.m} and Hypothesis~\ref{h.a}, where the Malthusian parameter associated with $\Pi^{(j)}$ is denoted by $p^*_j$. 
In addition, consider the $\mathcal S$-valued process  $(\lambda^I_\eta)_{\eta\in(0,1]}$ defined by
\[
\lambda^I_\eta:=\left(\lambda^I_{\eta,m}\right)_{m\in\N}:=\left(v_j\lambda^{(j)}_\eta\right)_{j\in\mathcal J},
\] 
where  $(\lambda^{(j)}_\eta)_{\eta\in(0,1]}$, with $\lambda^{(j)}_\eta=(\lambda^{(j)}_{\eta,k})_{k\in\N}$, denotes the terminal state of the fragmentation process  $\Pi^{(j)}$ stopped at the stopping line described in Section~\ref{s.hfp}.

\subsection{Main results}\label{s.mr}
The main purpose of the present paper is to lift the limit theorems described in Section~\ref{s.motivation} up to fragmentation processes with immigration.  To this end we shall resort to the set-up established in  Section~\ref{intro_immigration}. In order to state the results we  need the following condition regarding the total mass of immigrating particles:
\begin{equation}\label{e.fpwi.1b}
\sum_{j\in\mathcal J}v_j^{1+p^*_j}<\infty
\end{equation}
$\mathbb P$-almost surely. Note that this assumption is stronger than (\ref{e.immigration_measure}) for that it in particular implies that the total mass of all the immigrating blocks is finite. 

Let $f^{(j)}:\Omega\to\mathcal B^+$ for all $j\in\mathcal J$ and set $f^I:=(f^{(j)})_{j\in\mathcal J}$. In addition, let $\varphi:[0,1)\times\N\to\mathcal J$ be a random function that assigns to each $(\eta,m)\in(0,1]\times\N$ the index $j\in\mathcal J$ such that $\lambda^I_{\eta,m}$ corresponds to $\lambda^{(j)}_{\eta,k}$  for some $k\in\N$. In the setting with immigration the analogue of $\langle\rho_\eta,\cdot\rangle$, which was defined in (\ref{e.em.1}) in the context of standard fragmentation processes, is given as follows:
\[
\langle\rho^I_\eta,f^I\rangle := \sum_{m\in\N}\left(\lambda^I_{\eta,m}\right)^{1+p^*_{\varphi(\eta,m)}} f^{(\varphi(\eta,m))}\left(\frac{\lambda^I_{\eta,m}}{\eta}\right).
\]
Observe that
\begin{equation}\label{e.immigration}
\langle\rho^I_\eta,f^I\rangle = \sum_{j\in\mathcal J}\sum_{k\in\N}\left(v_j\lambda^{(j)}_{\frac{\eta}{v_j},k}\right)^{1+p^*_j}f^{(j)}\left(\frac{v_j\lambda^{(j)}_{\frac{\eta}{v_j},k}}{\eta}\right)= \sum_{j\in\mathcal J}v^{1+p^*_j}_j\langle\rho^{(j)}_\eta,f^{(j)}\rangle,
\end{equation}
where 
\begin{equation}\label{e.immigration.0c}
\langle\rho^{(j)}_\eta,f^{(j)}\rangle:=\sum_{k\in\N}\left(\lambda^{(j)}_{\frac{\eta}{v_j},k}\right)^{1+p^*_j}f^{(j)}\left(\frac{v_j\lambda^{(j)}_{\frac{\eta}{v_j},k}}{\eta}\right).
\end{equation}
for any $\eta\in(0,1]$. Notice that the definition (\ref{e.immigration.0c}) is  analogous  to (\ref{e.em.1}), but taking into account that the process starts from $v_j$. 
Furthermore, in view of (\ref{e.Lambda_0}) and (\ref{p.rho_f_b}) set
\[
\Lambda^I_0(f^I):=\sum_{j\in\mathcal J}v^{1+p^*_j}_j\langle\rho^{(j)},f^{(j)}\rangle\Lambda^{(j)}_0(p^*_j),
\]
where
\[
\langle\rho^{(j)},f^{(j)}\rangle=\lim_{\eta\downarrow0}\mathbb E\left(\langle\rho^{(j)}_\eta,f^{(j)}\rangle\right)\qquad\text{as well as}\qquad\Lambda^{(j)}_0(p^*_j):=\lim_{\eta\downarrow0}\langle\rho^{(j)}_\eta,{\bf 1}\rangle
\]
for all $j\in\mathcal J$. With this set-up we obtain the following result, whose proof  will be provided in Section~\ref{s.proofs}, which extends Proposition~\ref{t.1.1} to the setting with immigration. 

\begin{theorem}\label{t.1.2}
Let $f^I$ be as above. If (\ref{e.fpwi.1b}) as well as
\begin{equation}\label{e.L1_condition}
\sup_{j\in\mathcal J}\|f^{(j)}\|_\infty<\infty,
\end{equation} 
hold, then
\begin{equation}\label{e.t.SLLNfpi.1.0}
\lim_{\eta\downarrow0}\langle\rho^I_\eta,f^I\rangle=\Lambda^I_0(f^I)
\end{equation}
$\mathbb P$-almost surely. If in addition
\begin{equation}\label{e.fpwi.1b_L1}
\sum_{j\in\mathcal J}v_j^{1+p^*_j}\in\mathscr L^1(\mathbb P)
\end{equation}
is satisfied, then the convergence in (\ref{e.t.SLLNfpi.1.0}) also holds in
$\mathscr L^1(\mathbb P)$.
\end{theorem}
In view of Proposition~\ref{t.1.1} this theorem says that the limit of the series in (\ref{e.immigration}) as $\eta\downarrow0$ is the same as taking the limit inside the series. Note that this is not an obvious result, since in general neither the DCT nor the MCT is applicable in this situation. 

For the second main result of this paper  let $\phi^{(j)}:\rpn\times\mathcal P\times\Omega\to\R$, $j\in\mathcal J$, be  a random function such that $(\phi^{(j)}(x,\pi))_{x\in\rpn}$, $\pi\in\mathcal P$, is an $\R$-valued  stochastic processes,  with $\phi^{(j)}(x,\pi)=\phi^{(j)}(0,\pi)=0$ $\mathbb P$-a.s. for all $x>1$, which 
has c\`adl\`ag paths $\mathbb P$-a.s. and is independent of $\Pi$. In addition, we assume that $\phi^{(j)}(\cdot,(1,0,\ldots))=0$. 

The process $(Z^{I}_\eta)_{\eta\in(0,1]}$ counted with the random characteristics  $(\phi^{(j)})_{j\in\mathcal J}$ is then given by
\[
Z^{I}_\eta=\sum_{j\in\mathcal J}\sum_{i\in\mathcal I^{(j)}}\phi^{(j)}_i\left(\frac{\eta}{|\Pi^{(j)}_{\kappa^{(j)}(t_{j,i})}(t_{j,i}-)|},\pi^{(j)}(t_{j,i})\right)\mathds1_{\left\{|\Pi^{(j)}_{\kappa^{(j)}(t_{j,i})}(t_{j,i}-)|\ne0\right\}}
\]
for every $\eta\in(0,1]$, where the $\phi^{(j)}_i$ are independent copies of $\phi^{(j)}$ and where $(\kappa^{(j)}(t_{j,i}),\pi^{(j)}(t_{j,i}))_{i\in\mathcal I^{(j)}}$ are the atoms of the Poisson point process underlying $\Pi^{(j)}$. For the second main result of this paper we assume that the $\Pi^{(j)}$ and $\phi^{(j)}$ are equal in distribution to some fragmentation process $\Pi$ and some random characteristic $\phi$ respectively. Hence, we omit the index $j$ for terms that are determined by the distribution of $\Pi^{(j)}$ or $\phi^{(j)}$.
\begin{theorem}\label{t.rci}
Assume that the $(\Pi^{(j)})_{j\in\mathcal J}$ and $(\phi^{(j)})_{j\in\mathcal J}$ are each \iid and satisfy (\ref{e.assumption_t.rc}) and (\ref{e.assumption_1}),
respectively.
Moreover,  assume that (\ref{e.fpwi.1b_L1}) holds and that $\Pi^{(j)}$, $j\in\mathcal J$, is independent of  $(v_j)_{j\in\mathcal J}$ and $(\mathfrak t_j)_{j\in\mathcal J}$. Then
\begin{equation}\label{e.a.s.convergence}
\eta^{1+p^*}Z^I_\eta\to\sum_{j\in\mathcal J}v_j^{1+p^*} \frac{\Lambda^{(j)}_0(p^*)}{\Phi'(p^*)}\int_{(0,1)}\mathbb E\left(\sum_{k\in\N}|\Pi_{k}(t)|^{1+p^*}\int_\mathcal P\int_{(0,1]}s^{p^*}\phi\left(\rho,\pi\right)\dd s\mu(\dd\pi)\right)\dd t
\end{equation}
$\mathbb P$-a.s. and in $\mathscr L^1(\mathbb P)$ as $\eta\downarrow0$.
\end{theorem}

As an example where the above result can be applied we obtain the following corollary which can be proved along the lines of the proof of Theorem~3 of \cite{Kno12}, resorting to Theorem~\ref{t.rci}.

\begin{corollary}\label{c.fei}
Based on Theorem~\ref{t.rci} the statement of Theorem~3 of \cite{Kno12}, which is concerned with the asymptotic behaviour of the energy cost in the fragmentation energy model introduced in \cite{BM05}, can be extended to the setting with immigration. 
\end{corollary}

As mentioned before, the above results can be extended from homogeneous to self-similar fragmentation processes without any further effort.

\begin{corollary}\label{c.self.sim}
According to the same arguments which we used to extend Theorem~1 of \cite{HKK10} and Theorem~1 of \cite{Kno12} from the homogeneous to the self-similar case, also the statements of Theorem~\ref{t.1.2} and Theorem~\ref{t.rci} (as well as Corollary~\ref{c.fei})   remain true for self-similar fragmentation processes with indices of self-similarity in $\R$. Note that the reason for the independence of the abovementioned results from the indices of self-similarity is the fact that   the problems under consideration do not take into account the times at which fragmentation events occur. 
\end{corollary}

The final main result of this paper is concerned with the asymptotic decay of the size of the largest block in the process with immigration. More precisely, the following result gives bounds on the exponential decay rate. Let $\lambda^I_1(t)$ denote the size of the largest block present in the fragmentation process with immigration at time $t\in\rpn$. That is to say,
\[
\lambda^I_1(t):=
\max_{j\in\mathcal J}\left(v_j\lambda_1^{(j)}(t-\mathfrak t_j)\right).
\] 
Notice that the above maximum exists, since the total mass of all blocks at time $t\in\rpn$ is finite. 
Let us further emphasise  that in contrast to the above results the next result takes into account the time at which immigration events happen. 
For every $t\in\rpn$ set
\begin{equation}\label{e.j_t}
j_t:=\min\left\{j\in\mathcal J:\lambda^{(j)}_1(t-\mathfrak t_j)=\lambda^I_1(t)\right\}.
\end{equation}
Observe that this minimum exists, since all but finitely many blocks at time $t\in\rpn$ have size less than $\lambda^I_1(t)$. The following theorem, which we prove in Section~\ref{s.decay}, establishes bounds for the asymptotic decay rate of $\lambda^I_1(t)$ as $t\to\infty$.

\begin{theorem}\label{t.1.2b}
Assume that the $\Pi^{(j)}$, $j\in\mathcal J$, are identically distributed and that (\ref{e.fpwi.1b}) holds. 
Then we have
\[
-\Phi'(\bar p)\le\liminf_{t\to\infty}\frac{-\ln(\lambda^I_1(t))}{t}\le\limsup_{t\to\infty}\frac{-\ln(\lambda^I_1(t))}{t}\le-\Phi'(\bar p)+\limsup_{t\to\infty}\frac{\mathfrak t_{j_t}}{t}
\]
$\mathbb P$-almost surely.
\end{theorem}
Notice that $\nicefrac{\mathfrak t_{j_t}}{t}\le1$. Note further that if  $\lim_{t\to\infty}\nicefrac{\mathfrak t_{j_t}}{t}=0$, i.e. in particular if immigration is omitted and the fragmentation process starts from a random $\mathcal S$-valued vector ${\bf u}$, then in view of Corollary~1.4 of \cite{Ber06} 
the above theorem shows that the asymptotic decay of the largest block is the same as for fragmentation processes starting from a single block. 

\section{Example - spine decomposition}\label{s.example}
The aim of this section is to consider an example of a homogeneous  fragmentation process for which we can give an alternative proof that (\ref{e.t.SLLNfpi.1.0}) holds. This example is based on the spine decomposition for fragmentation processes.


Throughout this section fix some $p\in(\underline p,\infty)$ and recall that $\Pi$ is a homogeneous $\mathcal P$-fragmentation processes, which satisfies Hypothesis~\ref{h.m} and Hypothesis~\ref{h.a}, with dislocation measure $\mu$. Recall further that $\nu$ is the 
dislocation measure on $\mathcal S_1$ 
and that the measure $\mu$ on $\mathcal P$ is given by (\ref{e.mu}). In addition, let $p\in(\underline p,\infty)$ and consider the measure $\mu^{(p)}$ on $\mathcal P$ given by $\mu^{(p)}(d\pi) = |\pi_1|^{p}\mu(d\pi)$. Let $(\pi^1(t))_{t\in\rpn}$ be a Poisson point process on $\mathcal P$ with characteristic measure $\mu^{(p)}$ and let $(t_i)_{i\in\mathcal I_1}$ be the times for which this process takes a value in $\mathcal P\setminus\{(\N,\emptyset,\ldots)\}$, where $\mathcal I_1$ is an at most countable index set. Furthermore, let $\Pi^{(p)}$ be a standard homogeneous $\mathcal P$-fragmentation process under $\mathbb P$ with dislocation measure $\mu^{(p)}$ and such that the Poisson point process on $\mathcal P$ underlying $(\Pi^{(p)}_1(t))_{t\in\rpn}$ coincides with $(\pi^1(t))_{t\in\rpn}$. In addition, set
\[
\Delta(t):=\left|\Pi^{(p)}_1(t-)\right|\left|(\pi^1_n(t))_{n\in\N\setminus\{1\}}\right|^\downarrow
\]
for any $t\in\rpn$. Notice that $(\Delta(t))_{t\in\rpn}$ is a Poisson point process on $\mathcal S_1$ whose atoms in $\mathcal S_1\setminus\{(0,\ldots)\}$ are given by $(\Delta(t_i))_{i\in\mathcal I_1}$. Let $\lambda^{\Delta(t_i)}$, $i\in\mathcal I_1$, be  independent fragmentation processes, each starting from $\Delta(t_i)$, with dislocation measure $\nu$. Consider the fragmentation process with immigration $\lambda^I:=(\lambda^I(t))_{t\in\rpn}$ defined by
\begin{equation}\label{e.examplefragimmigration}
\lambda^I(t)=\sum_{i\in\mathcal I_1:t_i\le t}\lambda^{\Delta(t_i)}(t-t_i)
\end{equation}
for all $t\in\rpn$. Observe that this process starts from ${\bf u}=(0,\ldots)$. 

Recall that $\mathscr F$ is the filtration generated by $\Pi$ and for any $t\in\rpn$ consider the change of measure 
\begin{equation}\label{change_measure}
\left.\frac{\dd\mathbb P^{(p)}}{\dd\mathbb P}\right|_{\mathscr F_t}=e^{\Phi(p)t-p\xi(t)}.
\end{equation}
\begin{rem}\label{r.equivalentmeasures}
Assume that $p\in(\underline p,\bar p)$ and recall the martingale $M(p)$ given by (\ref{e.martingale_M}). We remark that considering projections of the  change of measure in (\ref{change_measure}) onto the sub-filtration $\mathscr G$, which is generated by the asymptotic frequencies of $\Pi$, results in
\[
\left.\frac{\dd\mathbb P^{(p)}}{\dd\mathbb P}\right|_{\mathscr G_t}
=M_t(p)
\]
for every $t\in\rpn$.
According to Theorem~1 of \cite{BR03} (cf. also Theorem~4 of \cite{BR05} for the conservative case) the unit-mean martingale $M(p)$ is uniformly integrable. Hence, $\mathbb E(M_\infty(p))=1$ and thus $\mathbb P^{(p)}$ is a probability measure on $\mathcal G_\infty:=\bigcup_{t\in\rpn}\mathcal G_t$.
Moreover, using that $\mathbb E(M_\infty(p))>0$ one obtains that $M_\infty(p)>0$ $\mathbb P$-a.s., see Lemma~1.35 of \cite{Kno11} (or Theorem 2 of  \cite{Ber03} for the conservative case). Consequently, restricted to the $\sigma$-algebra $\mathcal G_\infty$, the measures $\mathbb P^{(p)}$ and $\mathbb P$ are equivalent. 
\end{rem}

\begin{proposition}\label{l.e.fpi.1}
Let $f\in\mathcal B^+_b$. Then the process $\lambda^I$ constructed in (\ref{e.examplefragimmigration}) satisfies 
\[
\lim_{\eta\downarrow0}\langle\rho^I_\eta,f^I\rangle=\Lambda^I_0(f^I)
\]
$\mathbb P$-a.s., where $f^{(j)}:\equiv f$ for all $j\in\mathcal J$.
\end{proposition}

\begin{pf}
In view of (\ref{e.lambda}) let $\lambda$ be the mass fragmentation associated with $\Pi$ process. Furthermore, consider the following spine decomposition: 
\[
|\Pi(t)|=(|\Pi_1(t)|,0,\ldots)+\sum_{i\in\mathcal I_1:t_i\le t}\sum_{j\in\N\setminus\{1\}}\left|\Pi^{(i,j)}(t-t_i)\right|
\]
$\mathbb P^{(p)}$-a.s., where the $\Pi^{(i,j)}$ are independent and satisfy
\[
\mathbb P^{(p)}\left(\left.|\Pi^{(i,j)}(u-t_i)|\in\cdot\,\right|\mathscr F^1_{t_i}\right)=\left.\mathbb P^{(p)}\left(x_{i,j}|\Pi(u-t_i)|\in\cdot\right)\right|_{x_{i,j}=|\Pi_1(t_i-)\cap\pi_j(t_i)|}
\]
$\mathbb P^{(p)}$-a.s., where $(\mathscr F^1_t)_{t\in\rpn}$ is the filtration generated by $\Pi_1$. Moreover, under $\mathbb P^{(p)}$ the behaviour of the block $\Pi_1$, which is considered to be the {\it spine} or {\it tagged fragment}, is determined by a Poisson point process with intensity $\mu^{(p)}$.

Recall the construction of $\lambda^I$ in (\ref{e.examplefragimmigration}) and observe that 
\[
\lambda(t)=|\Pi_1(t)|+\lambda^I(t)
\]
$\mathbb P^{(p)}$-almost surely. 
That is to say, under $\mathbb P^{(p)}$ we can interpret the immigrating particles of $\lambda^I(t)$ as those particles that result from the fragmentation of the spine $\Pi_1$ at the jump times $(t_i)_{i\in\mathcal I_1}$ except for the tagged fragments $\Pi_1(t_i)$, $i\in\mathcal I_1$, themselves.

Recall the definition of $\upsilon_{\eta,1}$ in (\ref{e.stopping_line}). Using notations introduced in Section~\ref{intro_immigration} we infer from Proposition~\ref{t.1.1} that
\begin{align}\label{e.cI.0}
\lim_{\eta\downarrow0}\langle\rho^I_\eta,f^I\rangle &= \lim_{\eta\downarrow0}\langle\rho_\eta,f\rangle-\lim_{\eta\downarrow0}\left[|\Pi_1(\upsilon_{\eta,1})|^{1+p^*}f\left(\frac{|\Pi_1(\upsilon_{\eta,1})|}{\eta}\right)\right]\notag
\\[0.5ex]
&= \lim_{\eta\downarrow0}\langle\rho_\eta,f\rangle\notag
\\[0.5ex]
&= \langle\rho,f\rangle\Lambda(p^*)
\\[0.5ex]
&= \Lambda^I_0(f^I)\notag
\end{align}
$\mathbb P^{(p)}$-almost surely. Note that in order to apply Proposition~\ref{t.1.1} we have used that $\mathbb P^{(p)}$ and $\mathbb P$ are equivalent measures on $\mathscr G_\infty$, cf. Remark~\ref{r.equivalentmeasures}, to deduce that the convergence in Proposition~\ref{t.1.1} holds true $\mathbb P^{(p)}$-almost surely. Since the event $\{\lim_{\eta\downarrow0}\langle\rho^I_\eta,f^I\rangle=\Lambda^I_0(f^I)\}$ is $\mathscr G_\infty$-measurable, we conclude in view of (\ref{e.cI.0}) and resorting again to the fact that $\mathbb P^{(p)}$ and $\mathbb P$ are equivalent measures on $\mathscr G_\infty$ that 
\[
\lim_{\eta\downarrow0}\langle\rho^I_\eta,f^I\rangle=\Lambda^I_0(f^I)
\]
$\mathbb P$-almost surely.
\end{pf}

We remark that it follows from Lemma~2 in \cite{BR03} that (\ref{e.fpwi.1b}) is satisfied for the process $\lambda^I$  given by (\ref{e.examplefragimmigration}) and thus the statement of Proposition~\ref{l.e.fpi.1} also follows from Theorem~\ref{t.1.2}.

Let us now assume that the dislocation measure $\nu$ is conservative and let us finish this section by having a closer look at the characteristic measure under $\mathbb P^{(p)}$ of the Poisson measure $N$ that describes the immigration structure of (\ref{e.examplefragimmigration}). Note that $N$ is a random measure on $\mathcal S_1\otimes\rpn$ with atoms $(|(\pi_n(t_i))_{n\in\N\setminus\{1\}}|^\downarrow)_{i\in\mathcal I_1}$ in $\mathcal S_1\setminus\{(0,\ldots)\}$. The first thing to mention is that under $\mathbb P^{(p)}$ the intensity of $N$ is of the form $\nu_I\otimes\dd t$, where $\nu_I$ is a $\sigma$-finite measure on $\mathcal S_1$ and $\dd t$ denotes the Lebesgue measure on $\rpn$. Further, recall that the Poisson point process on $\mathcal P$ with atoms $(\pi(t_i))_{i\in\mathcal I_1}$ in $\mathcal P\setminus(\N,\emptyset,\ldots)$ has characteristic measure $\mu^{(p)}$. Hence, since $\nu$ is conservative, the measure $\nu_I$ is the projection of $\mu^{(p)}$ on $\mathcal S_1$ and in view of (3) in \cite{HKK10} we thus infer that
\[
\int_{\mathcal S_1}g({\bf s})\nu_I(\dd{\bf s})=\int_\mathcal Pg(|\pi|^\downarrow)\mu^{(p)}(\dd\pi)=\int_{\mathcal P}g(|\pi|^\downarrow)|\pi_1|^p\mu(\dd\pi)=\int_{\mathcal S_1}g({\bf s})\sum_{n\in\N}s_n^{1+p}\nu(\dd{\bf s})
\] 
holds for any nonnegative test function $g:\mathcal S_1\to\rpn$, which results in
\[
\nu_I(\dd{\bf s})=\sum_{n\in\N}s_n^{1+p} \nu(\dd{\bf s})
\]
for all ${\bf s}\in\mathcal S_1$.

\section{Proof of Theorem~\ref{t.1.2}}\label{s.proofs}

The goal of this section is to prove Theorem~\ref{t.1.2}. In order to tackle the proof of this result we first need to develop some auxiliary lemmas. Let us mention that we make use of ideas of \cite{99} and \cite{Olo96}. 

Recall the 
 set-up established in Section~\ref{intro_immigration}. 
Throughout this section, unless stated otherwise, let $\mathcal J$ be a deterministic at most countably infinite index set and let $(v_j)_{j\in\mathcal J}\in\mathcal S$ be a deterministic vector. Note that this is somewhat an abuse of notation, since before we denoted by $\mathcal J$ and $(v_j)_{j\in\mathcal J}$ the random objects associated with the immigration process. However, our goal in this section is to prove the desired convergence for the abovementioned deterministic objects and, since we are interested in almost sure convergence, this then proves the result with the random objects we are interested in.

For every $\eta\in(0,1]$ set 
\begin{equation}\label{e.J_t}
\mathcal J_\eta:=\{j\in\mathcal J:v_j\ge \eta\}\qquad\text{as well as}\qquad\mathcal J^\complement_\eta:=\{j\in\mathcal J:v_j<\eta\}.
\end{equation}

For each $j\in\mathcal J$ and $\eta\in(0,1]$ let $(\lambda^{(j)}_\eta)_{\eta\in(0,1]}$ be the stopped fragmentation process associated with $\Pi^{(j)}$. In addition, let $(\mathscr H^{(j)}_\eta)_{\eta\in(0,1]}$\label{p.H-j}, $j\in\mathcal J$, be the filtration generated by the stopped process $(\lambda^{(j)}_\eta)_{\eta\in(0,1]}$, i.e. 
\[
\mathscr H^{(j)}_\eta=\sigma\left(\left\{\lambda^{(j)}_u:u\in\left[\nicefrac{\eta}{v_j},1\right]\right\}\right). 
\]
Notice that $v_j$ is $\mathscr H^{(j)}_\eta$-measurable for all $\eta\in(0,1]$ and consider the filtration $\mathscr H^I:=(\mathscr H^I_\eta)_{\eta\in(0,1]}$ given by 
\[
\mathscr H^I_\eta:=\sigma\left(\bigcup_{j\in\mathcal J_\eta}\mathscr H^{(j)}_\eta\right)
\] 
for any $\eta\in(0,1]$. 


The first lemma in this section in particular estalishes a useful submartingale property.

\begin{lemma}\label{l.fpwi.2}
Assume that (\ref{e.fpwi.1b})  holds.  Then there exists a $\Lambda^I_0\in\mathscr L^1(\mathbb P)$\label{p.Lambda^I_infty} such that $\langle\rho^I_\eta,{\bf 1}\rangle\to\Lambda^I_0$ $\mathbb P$-a.s. as $\eta\downarrow0$.
\end{lemma}

\begin{pf}
The idea of the proof is to use the submartingale convergence theorem.

By means of the MCT we infer from (\ref{e.fpwi.1b}) that
\begin{equation}\label{e.l.fpwi.2.1}
\sup_{\eta\in(0,1]}\mathbb E(\langle\rho^I_\eta,{\bf1}\rangle)=\sup_{\eta\in(0,1]}\sum_{j\in\mathcal J}v_j^{1+p^*_j}\mathbb E\left(\langle\rho^{(j)}_\eta,{\bf1}\rangle\right)=\sum_{j\in\mathcal J}v_j^{1+p^*_j}<\infty.
\end{equation}
Moreover, the MCT for conditional expectations in conjunction with the martingale property of $\Lambda^{(j)}(p^*)=(\langle\rho^{(j)}_\eta,{\bf 1}\rangle)_{\eta\in(0,1]}$, $j\in\mathcal J$, yields that
\begin{align*}
\mathbb E\left(\left.\left\langle\rho^I_{\eta\rho},{\bf 1}\right\rangle\right|\mathscr H^I_\eta\right) &= \sum_{j\in\mathcal J}v_j^{1+p^*_j}\mathbb E\left(\left.\left\langle\rho^{(j)}_{\eta\rho},{\bf 1}\right\rangle\right|\mathscr H^{(j)}_\eta\right)
\\[0.5ex]
&\ge \sum_{j\in\mathcal J_\eta}v_j^{1+p^*_j}\mathbb E\left(\left.\left\langle\rho^{(j)}_{\eta\rho},{\bf 1}\right\rangle\right|\mathscr H^{(j)}_\eta\right)
\\[0.5ex]
&= \sum_{j\in\mathcal J_\eta}v_j^{1+p^*_j}\left\langle\rho^{(j)}_\eta,{\bf 1}\right\rangle
\\[0.5ex]
&= \langle\rho^I_\eta,{\bf 1}\rangle
\end{align*}
$\mathbb P$-a.s. for all $\eta,\rho\in(0,1]$, which shows that under $\mathbb P$ the process $(\langle\rho^I_\eta,{\bf 1}\rangle)_{\eta\in(0,1]}$ is a nonnegative $\mathscr H^I$-submartingale. Note that here we have used the independence of $(\Pi^{(j)})_{j\in\mathcal J}$. In view of (\ref{e.l.fpwi.2.1}) the submartingale convergence theorem thus ensures that  there exists a $\Lambda^I_0\in\mathscr L^1(\mathbb P)$ such that $\langle\rho^I_\eta,{\bf 1}\rangle\to\Lambda^I_0$ $\mathbb P$-a.s. as $\eta\downarrow0$. 
\end{pf}

The previous lemma can be strengthened in the sense that the obtained limiting random variable can be described explicitly. Indeed, this assertion is the statement of the following proposition. 

\begin{proposition}\label{p.SLLNfpi.2}
Assume that (\ref{e.fpwi.1b})  holds. Then we have 
\[
\langle\rho^I_\eta,{\bf 1}\rangle\to\sum_{j\in\mathcal J}v_j^{1+p^*_j}\Lambda^{(j)}_0(p^*_j)
\] 
$\mathbb P$-a.s. as $\eta\downarrow0$.
\end{proposition}  

\begin{pf}
As a consequence of the MCT and (\ref{e.fpwi.1b}) we obtain analogously to (\ref{e.l.fpwi.2.1}) that
\begin{equation}\label{e.fpwi.1}
\mathbb E\left(\sum_{j\in\mathcal J}v_j^{1+p^*_j}\Lambda^{(j)}_0(p^*_j)\right)=\sum_{j\in\mathcal J}v_j^{1+p^*_j}\mathbb E\left(\Lambda^{(j)}_0(p^*_j)\right)=\sum_{j\in\mathcal J}v_j^{1+p^*_j}<\infty,
\end{equation}
where we have used that the unit-mean martingale $\Lambda^{(j)}(p^*_j)$ is uniformly integrable. Further, let $\Lambda^I_0$ be given by Lemma~\ref{l.fpwi.2} and recall the definition of $\mathcal J_\eta$ as well as $\mathcal J^\complement_\eta$ in (\ref{e.J_t}). 
Observe that for any $0<\eta\le\rho\le1$ we have 
\begin{align}\label{e.p.SLLNfpi.2.1}
&\Lambda^I_0-\sum_{j\in\mathcal J}v_j^{1+p^*_j}\Lambda^{(j)}_0(p^*_j)
\\[0.5ex]
&= \Lambda^I_0-\langle\rho^I_\eta,{\bf 1}\rangle+
\sum_{j\in\mathcal J_\rho}v_j^{1+p^*_j}\left(\left\langle\rho^{(j)}_\eta,{\bf1}\right\rangle-\Lambda^{(j)}_0(p^*_j)\right)
 +\sum_{j\in\mathcal J^\complement_\rho}v_j^{1+p^*_j}\left\langle\rho^{(j)}_\eta,{\bf1}\right\rangle-\sum_{j\in\mathcal J^\complement_\rho}v_j^{1+p^*_j}\Lambda^{(j)}_0(p^*_j).\notag
\end{align}
According to Lemma~\ref{l.fpwi.2} we have that
\begin{equation}\label{e.fpwi.2}
\Lambda^I_0-\langle\rho^I_\eta,{\bf 1}\rangle\to0
\end{equation}
as $\eta\downarrow0$. Note further that by means of  Proposition~\ref{t.1.1} and (\ref{e.fpwi.1b}) we have 
\begin{equation}\label{e.fpwi.2_a}
\lim_{\eta\downarrow0}\sum_{j\in\mathcal J_\rho}v_j^{1+p^*_j}\left(\left\langle\rho^{(j)}_\eta,f^{(j)}\right\rangle-\left\langle\rho^{(j)},f^{(j)}\right\rangle\Lambda^{(j)}_0(p^*_j)\right)=0
\end{equation}
$\mathbb P$-a.s. for any $\rho\in(0,1]$, where we have used that this sum has only finitely many summands as infinitely many $j\in\mathcal J$ with $v_j\ge \rho$ would contradict (\ref{e.fpwi.1b}). Moreover, resorting to (\ref{e.fpwi.1}) we obtain that 
\begin{equation}\label{e.fpwi.2_b}
\sum_{j\in\mathcal J^\complement_\rho}v_j^{1+p^*_j}\Lambda^{(j)}_0(p^*_j)\to0
\end{equation}
as $\rho\downarrow0$. Let us now consider the penultimate term in (\ref{e.p.SLLNfpi.2.1}). To this end, notice that the limit $\lim_{\eta\downarrow0}\sum_{j\in\mathcal J^\complement_\rho}v_j^{1+p^*_j}\langle\rho^{(j)}_\eta,{\bf1}\rangle$ exists $\mathbb P$-a.s., since, according to Lemma~\ref{l.fpwi.2}, $\lim_{\eta\downarrow0}\langle\rho^{I}_\eta,{\bf1}\rangle$ exists $\mathbb P$-a.s. and also $\lim_{\eta\downarrow0}\sum_{j\in\mathcal J_\rho}v_j^{1+p^*_j}\langle\rho^{(j)}_\eta,{\bf1}\rangle$ exists, because the sum is taken over only finitely many summands. Since the map
\[
\rho\mapsto\lim_{\eta\downarrow0}\sum_{j\in\mathcal J^\complement_\rho}v_j^{1+p^*_j}\langle\rho^{(j)}_\eta,{\bf1}\rangle
\]
is monotone in $\rho$, we infer that the limit as $\rho\downarrow0$ exists $\mathbb P$-almost surely. Hence,  we deduce by means of Fatou's lemma and the martingale property of $\Lambda^{(j)}(p^*_j)=(\langle\rho^{(j)}_\eta,{\bf1}\rangle)_{\eta\in(0,1]}$ that
\[
\mathbb E\left(\lim_{\rho\downarrow0}\lim_{\eta\downarrow0}\sum_{j\in\mathcal J^\complement_\rho}v_j^{1+p^*_j}\langle\rho^{(j)}_\eta,{\bf1}\rangle\right) \le \liminf_{\rho\downarrow0}\liminf_{\eta\downarrow0}\sum_{j\in\mathcal J^\complement_\rho}v_j^{1+p^*_j}\mathbb E\left(\langle\rho^{(j)}_\eta,{\bf1}\rangle\right)
\le \lim_{\rho\downarrow0}\sum_{j\in\mathcal J^\complement_\rho}v_j^{1+p^*_j}=0
\]
$\mathbb P$-a.s., since $\sum_{j\in\mathcal J}v_j^{1+p^*_j}<\infty$. Consequently, as
\[
\lim_{\rho\downarrow0}\lim_{\eta\downarrow0}\sum_{j\in\mathcal J^\complement_\rho}v_j^{1+p^*_j}\langle\rho^{(j)}_\eta,{\bf1}\rangle\ge0,
\]
this implies that
\[
\lim_{\rho\downarrow0}\lim_{\eta\downarrow0}\sum_{j\in\mathcal J^\complement_\rho}v_j^{1+p^*_j}\langle\rho^{(j)}_\eta,{\bf1}\rangle=0
\]
$\mathbb P$-almost surely. Combining this with (\ref{e.fpwi.2}), (\ref{e.fpwi.2_a}) and (\ref{e.fpwi.2_b}) it thus follows from (\ref{e.p.SLLNfpi.2.1}) that
\[
\Lambda^I_0-\sum_{j\in\mathcal J}v_j^{1+p^*_j}\Lambda^{(j)}_0(p^*_j)=0
\]
$\mathbb P$--a.s., which completes the proof.
\end{pf}

We are now in a position to prove Theorem~\ref{t.1.2}. 

{\bf Proof of Theorem~\ref{t.1.2}}
We first prove almost sure convergence in (\ref{e.t.SLLNfpi.1.0}) and then also the corresponding $\mathscr L^1$-convergence follows along the lines of  these arguments.

\underline{Part I}
By means of  Proposition~\ref{t.1.1} and Fatou's lemma we have
\begin{equation}\label{e.t.SLLN.fpi.1}
\liminf_{\eta\downarrow0}\langle\rho^I_\eta,f^I\rangle\ge\sum_{j\in\mathcal J}v_j^{1+p^*_j}\langle\rho^{(j)},f^{(j)}\rangle\Lambda^{(j)}_0(p^*_j)
\end{equation}
$\mathbb P$-almost surely. In view of (\ref{e.L1_condition}) set $f^I_\infty:=(\sup_{i\in\mathcal J}\|f^{(i)}\|_\infty)_{j\in\mathcal J}$. As a consequence of the additivity of $\langle\rho^{(j)},\cdot\rangle$ and $\langle\rho^{(j)},{\bf 1}\rangle=1$ for all $j\in\mathcal J$, we infer from Proposition~\ref{p.SLLNfpi.2} and (\ref{e.t.SLLN.fpi.1}), applied to  $f^I_\infty-f^I$ rather than $f^I$, that 
\begin{align*}
\limsup_{\eta\downarrow0}\langle\rho^I_\eta,f^I\rangle
&\le \lim_{\eta\downarrow0}\langle\rho^I_\eta,f^I_\infty\rangle-\liminf_{\eta\downarrow0}\langle\rho^I_\eta,f^I_\infty-f^I\rangle
\\[0.5ex]
&\le \sum_{j\in\mathcal J}v_j^{1+p^*_j}\langle\rho^{(j)},f^I_\infty\rangle\Lambda^{(j)}_0(p^*_j)-\sum_{j\in\mathcal J}v_j^{1+p^*_j}\langle\rho^{(j)},f^I_\infty-f^{(j)}\rangle\Lambda^{(j)}_0(p^*_j)
\\[0.5ex]
&= \sum_{j\in\mathcal J}v_j^{1+p^*_j}\langle\rho^{(j)},f^{(j)}\rangle\Lambda^{(j)}_0(p^*_j)
\end{align*}
$\mathbb P$-a.s., which combined with (\ref{e.t.SLLN.fpi.1}) proves that
\begin{equation}\label{e.t.SLLN.fpi.2}
\lim_{\eta\downarrow0}\langle\rho^I_\eta,f^I\rangle=\sum_{j\in\mathcal J}v_j^{1+p^*_j}\langle\rho^{(j)},f^{(j)}\rangle\Lambda^{(j)}_0(p^*_j)
\end{equation}
$\mathbb P$-almost surely. Observe that this proves that the convergence in (\ref{e.t.SLLNfpi.1.0}) holds $\mathbb P$-a.s., since replacing the deterministic objects $\mathcal J$ and $(v_j)_{j\in\mathcal J}$ by the homonymous random objects of Section~\ref{intro_immigration} does not affect the almost sure convergence in (\ref{e.t.SLLN.fpi.2}).

Now assume that  $\mathcal J$ as well as $(v_j)_{j\in\mathcal J}$ are the random objects defined in Section~\ref{intro_immigration}. Following the lines of the above proof, resorting to 
\[
\lim_{\eta\downarrow0}\mathbb E\left(\langle\rho^I_\eta,{\bf1}\rangle\right)=\lim_{\eta\downarrow0}\sum_{j\in\N}\mathbb E\left(\mathds1_{\{j\in\mathcal J\}}v_j^{1+p^*_j}\right)\mathbb E\left(\langle\rho^{(j)}_\eta,{\bf1}\rangle\right)=\mathbb E\left(\sum_{j\in\mathcal J}v_j^{1+p^*_j}\right)
\]
instead of Proposition~\ref{p.SLLNfpi.2}, we infer that 
\[
\lim_{\eta\downarrow0}\mathbb E\left(\langle\rho^I_\eta,f^I\rangle\right)=\sum_{j\in\mathcal J}v_j^{1+p^*_j}\langle\rho^{(j)},f^{(j)}\rangle,
\]
which by means of (\ref{e.t.SLLN.fpi.2}) and Lemma~21.6 in \cite{Bau01} proves that (\ref{e.t.SLLNfpi.1.0}) also holds in $\mathscr L^1(\mathbb P)$.
\hfill$\square$

\begin{rem}\label{r.t.1.2.1}
Assume that the $\Pi^{(j)}$, $j\in\mathcal J$, are \iid and thus we omit the index $j$ for terms which are determined by the distribution of $\Pi^{(j)}$.
Then (\ref{e.t.SLLN.fpi.2}) can be proven without resorting to Lemma~\ref{l.fpwi.2} and Proposition~\ref{p.SLLNfpi.2}. Indeed, resorting to the DCT we obtain that
\begin{align}\label{e.r.t.1.2.1.3}
\lim_{\eta\downarrow0}\left\langle\rho^I_\eta,f^I\right\rangle &= \lim_{\eta\downarrow0}\sum_{j\in\mathcal J}v^{1+p^*}_j\langle\rho^{(j)}_\eta,f^{(j)}\rangle\notag
\\[0.5ex]
&= \sum_{j\in\mathcal J}v^{1+p^*}_j\lim_{\eta\downarrow0}\langle\rho^{(j)}_\eta,f^{(j)}\rangle
\\[0.5ex]
&= \sum_{j\in\mathcal J}v_j^{1+p^*_j}\langle\rho^{(j)},f^{(j)}\rangle\Lambda^{(j)}_0(p^*_j)\notag
\end{align}
holds $\mathbb P$-a.s., which proves  (\ref{e.t.SLLN.fpi.2}) in this special situation. Note that in (\ref{e.r.t.1.2.1.3}) we can indeed resort to the DCT, since an application of the MCT yields that
\[
\mathbb E\left(\sup_{\eta\in(0,1]}\sum_{j\in\mathcal J}v_j^{1+p^*}\langle\rho^{(j)}_\eta,f^{(j)}\rangle\right) \le \sup_{j\in\mathcal J}\|f^{(j)}\|_\infty\mathbb E\left(\sup_{\eta\in(0,1]}\langle\rho_\eta,{\bf 1}\rangle\right)\sum_{j\in\mathcal J}v^{1+p^*}_j< \infty,
\]
where the finiteness is a consequence of (\ref{e.fpwi.1b}), (\ref{e.L1_condition}) and Proposition~3.5 in \cite{Kno11}. 
\end{rem}

\section{Proof of Theorem~\ref{t.rci}}\label{s.proofs.2}
In this section we prove Theorem~\ref{t.rci}. For this purpose we work with the set-up established in  Section~\ref{intro_immigration}. As in Section~\ref{s.proofs} we consider throughout this section, unless stated otherwise, a deterministic vector $(v_j)_{j\in\mathcal J}\in\mathcal S$, which, as before, is sufficient even though Theorem~\ref{t.rci} is concerned with a random vector $(v_j)_{j\in\mathcal J}$.

In order to tackle the proof of Theorem~\ref{t.rci} we shall need the following result:

\begin{proposition}\label{l.t.SLLN.fprc.1b.new.1b.0}
Let $(\Pi^{(j)})_{j\in\mathcal J}$ and $(\phi^{(j)})_{j\in\mathcal J}$ be as in Theorem~\ref{t.rci}. Then we have
\[
\lim_{\eta\downarrow0}\sum_{j\in\mathcal J}v_j^{1+p^*}\eta^{1+p^*+\beta}Z^{\phi^{(j)}}_{\eta}=0
\]
$\mathbb P$-a.s. for all $\beta>0$.
\end{proposition}

\begin{pf}

Let $a,\beta>0$ as well as $\rho\in(0,1)$ and let $\phi$ be a generic random characteristics which in distribution equals $\phi^{(j)}$, $j\in\mathcal J$.  Furthermore, for every $j\in\mathcal J$ set
\[
Y^{\phi^{(j)}}_{a,\beta}:=\sum_{i\in\mathcal I}\mathds1_{\{t_i\le a\}}\sup_{\eta\in(0,1]}\left(\eta^{1+p^*+\beta}\phi^{(j)}_i\left(\eta,\pi(t_i)\right)\right),
\]
where the $\phi^{(j)}_i$  are independent copies of $\phi$. 
Observe that
\begin{equation}\label{Y_finite}
\mathbb E\left(Y^\phi_{a,\beta}\right)<\infty.
\end{equation}
Indeed, recall that  $\phi$ and $\Pi$ are independent.  Hence, the compensation formula for Poisson point processes and Tonelli's theorem yield that
\begin{align*}
& \mathbb E\left(\sum_{i\in\mathcal I}\mathds1_{\{t_i\le a\}}\sup_{\eta\in(0,1]}\left(\eta^{1+p^*+\beta}\phi_i\left(\eta,\pi(t_i)\right)\right)\right)
\\[0.5ex]
&= a\int_\mathcal P\mathbb E\left(\sup_{\eta\in(0,1]}\left(\eta^{1+p^*+\beta}\phi\left(\eta,\pi\right)\right)\right)\mu(\dd\pi)
\\[0.5ex]
&< \infty,
\end{align*}
where the $\phi_i$ are independent copies of $\phi$, which proves that (\ref{Y_finite}) holds.

In order to proceed with the proof we need to introduce some more notation. For this purpose, recall that the $(\kappa^{(j)}(t_{j,i}),\pi^{(j)}(t_{j,i}))_{i\in\mathcal I^{(j)}}$, $j\in\mathcal J$, are the atoms of the Poisson point process underlying $\Pi^{(j)}$. For any $j\in\mathcal J$ let $\mathcal I^{(j,a)}$ be given by
\[
\mathcal I^{(j,a)}:=\left\{i\in\mathcal I^{(j)}:t_{j,i}\in(na)_{n\in\N}\right\}
\]
and set 
\[
\mathcal L^{(j,a)}_{\eta,m}:=\left\{i\in\mathcal I^{(j,a)}:\left(|\Pi^{(j)}_{\kappa^{(j)}(t_{j,i})}(t_{j,i}-)|\ge \eta\right)\land\left(t_{j,i}= ma\right)\right\}
\]
for all $\eta\in(0,1]$ and $m\in\N$. By means of the strong fragmentation property we then have that
\begin{align*}
& \sum_{j\in\mathcal J}v^{1+p^*}_j\rho^{r(1+p^*+2\beta)}Z^{\phi^{(j)}}_{\rho^r}\notag
\\[0.5ex]
&\le \sum_{j\in\mathcal J}v^{1+p^*}_j\sum_{k=1}^{\lceil r\rceil}\sum_{m=1}^\infty\sum_{l\in\mathcal L^{(j,a)}_{\rho^{k},m}\setminus\mathcal L^{(j,a)}_{\rho^{k-1},m-1}}\sum_{i\in\mathcal I^{(j,l)}}\mathds1_{\{t_{j,i,l}\le a\}}\rho^{r(1+p^*+2\beta)}\phi^{(j)}_{i,l}\left(\frac{\rho^r}{|\Pi_{\kappa(t_{j,i,l})}(t_{j,i,l}-)|},\pi(t_{j,i,l})\right)\notag
\\[0.5ex]
&\le \rho^{r\beta}\rho^{-(1+p^*+\beta)}\sum_{j\in\mathcal J}v^{1+p^*}_j\sum_{k=1}^{\lceil r\rceil}\sum_{m=1}^k\sum_{l\in\mathcal L^{(j,a)}_{\rho^{k},m}\setminus\mathcal L^{(j,a)}_{\rho^{k-1},m-1}}\rho^{k(1+p^*+\beta)}Y^{(j,l)}_{a,\beta}
\end{align*}
$\mathbb P$--a.s. for each $r\in\rpn$, where 
the $(t_{j,i,l})_{i\in\mathcal I^{(j,l)}}$, $\phi^{(j)}_{i,l}$ and $Y^{(j,l)}_{a,\beta}$ are independent copies of $(t_i)_{i\in\mathcal I}$, $\phi$ and $Y^\phi_{a,\beta}$ respectively.  
Note that, by Tomelli's theorem,
\begin{align}\label{e.chain_inequalities}
& \mathbb E\left(\sum_{j\in\mathcal J}v^{1+p^*}_j\sum_{k=1}^{\infty}\sum_{m=1}^k\sum_{l\in\mathcal L^{(j,a)}_{\rho^{k},m}\setminus\mathcal L^{(j,a)}_{\rho^{k-1},m-1}}\rho^{k(1+p^*+\beta)}Y^{(j,l)}_{a,\beta}\right)\notag
\\[0.5ex]
&= \mathbb E\left(\sum_{j\in\mathcal J}v^{1+p^*}_j\sum_{k=1}^{\infty} \rho^{k(1+p^*+\beta)}\sum_{m=1}^k\sum_{l\in\N}\mathbb E\left(\mathds1_{l\in\mathcal L^{(j,a)}_{\rho^{k},m}\setminus\mathcal L^{(j,a)}_{\rho^{k-1},m-1}}Y^{(j,l)}_{a,\beta}\right)\right)\notag
\\[0.5ex]
&= \mathbb E\left(Y_{a,\beta}\right)\mathbb E\left(\sum_{j\in\mathcal J}v^{1+p^*}_j\sum_{k=1}^{\infty} \rho^{k(1+p^*+\beta)}\sum_{m=1}^k\mathbb E\left(\sharp\left(\mathcal L^{(j,a)}_{\rho^{k},m}\setminus\mathcal L^{(j,a)}_{\rho^{k-1},m-1}\right)\right)\right)
\\[0.5ex]
&= \mathbb E\left(Y_{a,\beta}\right)\mathbb E\left(\sum_{j\in\mathcal J}v^{1+p^*}_j\sum_{k=1}^{\infty} \rho^{k\beta}\mathbb E\left(\rho^{k(1+p^*)}T^{(a)}_{\rho^{k}}\right)\right)\notag
\\[0.5ex]
&\le \mathbb E\left(Y_{a,\beta}\right)\mathbb E\left(\sup_{u\in(0,1]}u^{1+p^*}T^{(a)}_u\right)\mathbb E\left(\sum_{j\in\mathcal J}v^{1+p^*}_j\right)\sum_{k=1}^{\infty} \rho^{k\beta},\notag
\end{align}
where 
\[
T^{(a)}_\eta:=\sharp\left\{i\in\mathcal I:\left(t_i\in(na)_{n\in\N}\left)\land\right(|\Pi_{\kappa(t_i)}(t_i-)|\ge\eta\right)\right\}
\]
for every $\eta\in(0,1]$.
Since the geometric series $\sum_{k=1}^{\infty} \rho^{k\beta}$ converges, we have according to (\ref{Y_finite}), Proposition~5 of \cite{Kno12} and  (\ref{e.fpwi.1b_L1})  that all the factors on the right-hand side of (\ref{e.chain_inequalities}) are finite. Consequently, we conclude that
\begin{align*}
& \mathbb E\left(\limsup_{\eta\downarrow0}\sum_{j\in\mathcal J}v_j^{1+p^*}\eta^{1+p^*+\beta}Z^{\phi^{(j)}}_{\eta}\right)
\\[0.5ex]
&= \mathbb E\left(\limsup_{r\to\infty}\sum_{j\in\mathcal J}v^{1+p^*}_j\rho^{r(1+p^*+2\beta)}Z^{\phi^{(j)}}_{\rho^r}\right) 
\\[0.5ex]
&\le \lim_{r\to\infty} \rho^{r\beta}\rho^{-(1+p^*+\beta)}\mathbb E\left(\sum_{j\in\mathcal J}v^{1+p^*}_j\sum_{k=1}^{\infty}\sum_{m=1}^k\sum_{l\in\mathcal L^{(j,a)}_{\rho^{k},m}\setminus\mathcal L^{(j,a)}_{\rho^{k},m-1}}\rho^{k\beta}Y^{(j,l)}_{a,\beta}\right)
\\[0.5ex]
&= 0,
\end{align*}
which proves the assertion, since $\sum_{j\in\mathcal J}v_j^{1+p^*}\eta^{1+p^*+\beta}Z^{\phi^{(j)}}_{\eta}\ge0$ for all $\eta\in(0,1]$.
\end{pf}

In the light of Proposition~\ref{p.SLLNfpi.2}  we also have the following extension of Lemma~3 in \cite{HKK10}:
\begin{lemma}\label{l.L_3}
Assume that (\ref{e.fpwi.1b_L1}) holds. 
Then there exists some $s_0\in(0,\infty)$ such that 
\[
\lim_{\eta\downarrow0}\sum_{j\in\mathcal J}v_j^{1+p^*}\sum_{k\in\N:\,\lambda^{(j)}_{\eta,k}<\eta^s}\left(\lambda^{(j)}_{\eta,k}\right)^{1+p^*}=0
\]
$\mathbb P$-a.s. for all $s\ge s_0$.
\end{lemma}

\begin{pf}
Let $p\in(\underline p,p^*)$ and observe  that the method of proof in Lemma~1 of \cite{HKK10} yields that also the  process given by
\[
\sum_{k\in\N}\left(\lambda^{(j)}_{\eta,k}\right)^{1+p}e^{\Phi^{(j)}(p)\sigma^{(j)}_{\eta,k}}
\]
for every $j\in\mathcal J$ and $\eta\in(0,1]$, where $\sigma^{(j)}_{\eta,k}$ denote the time of creation of  the block $\lambda^{(j)}_{\eta,k}$, is a unit-mean martingale. Hence, the arguments in the proof of Proposition~\ref{p.SLLNfpi.2} thus show that the process defined by
\[
\sum_{j\in\mathcal J}v_j^{1+p^*}\sum_{k\in\N}\left(\lambda^{(j)}_{\eta,k}\right)^{1+p}e^{\Phi^{(j)}(p)\sigma^{(j)}_{\eta,k}}
\]
for all $\eta\in(0,1]$ is a submartingale. 
Consequently, resorting to the submartingale convergence theorem and (\ref{e.fpwi.1b_L1}), which results in
\[
\sup_{\eta\in(0,1]}\mathbb E\left(\sum_{j\in\mathcal J}v_j^{1+p^*}\sum_{k\in\N}\left(\lambda^{(j)}_{\eta,k}\right)^{1+p}e^{\Phi^{(j)}(p)\sigma^{(j)}_{\eta,k}}\right)\le\sup_{\eta\in(0,1]}\mathbb E\left(\sum_{j\in\mathcal J}v_j^{1+p^*}\right)<\infty,
\] 
the lemma can be proven along the lines of Lemma~3 in \cite{HKK10}.
\end{pf}

{\bf Proof of Theorem~\ref{t.rci}}
Examining the proof of Theorem~1 in  \cite{Kno12} one realises that most proofs make use of the expectation of $\eta^{1+p^*}Z^\phi_\eta$. In view of (\ref{e.fpwi.1b_L1}) and the independence between the Poisson point processes governing the immigration and the fragmentation respectively, such arguments can be used also to deal with the process with immigration. Note that for this purpose we make use of the assumption that the immigrating fragmentation processes as well as the associated random characteristics are identically distributed, which results in  the expected values being the same for all $j\in\mathcal J$. However,  there are some steps where almost sure estimates are used and we need to use different arguments for those situations. Nonetheless, in the spirit of the above remarks, the proof of Theorem~\ref{t.rci} basically consists of showing that resorting to Proposition~\ref{l.t.SLLN.fprc.1b.new.1b.0} and Lemma~\ref{l.L_3} the method of proof used in \cite{Kno12} can be extended to  prove Theorem~\ref{t.rci}. In order to avoid replications of the arguments of \cite{Kno12} we merely outline how the proof of Theorem~1 of \cite{Kno12} can be adapted to the setting with immigration.

As before, for objects which are determined by the distribution of $\Pi^{(j)}$ or $\phi^{(j)}$ we omit the index $j\in\mathcal J$. For every $\iota>1$ and $\eta\in(0,1]$ and $j\in\mathcal J$ define a function $\phi^{(j)}_{\iota,\eta}:\rpn\times\mathcal P\times\Omega\to\R$
\[
\phi^{(j)}_{\iota,\eta}(x,\pi):=\phi^{(j)}(x,\pi)\mathds 1_{\{x>\eta^{\iota-1}\}}
\]
for all $x\in[0,1]$ and $\pi\in\mathcal P$. 
Let us first show that the analogue of Proposition~7 of \cite{Kno12} holds in the setting with immigration. That is to say,
\begin{equation}\label{e.Prop_8}
\lim_{k\to\infty}\sum_{j\in\mathcal J}v^{1+p^*}_j\rho^{k\delta\iota(1+p^*)}Z^{\phi^{(j)}_{\iota,\rho^{k\delta}}}_{\rho^{k\delta\iota}}=\lim_{\eta\downarrow0}\mathbb E\left(\eta^{1+p^*}Z^\phi_\eta\right)\sum_{j\in\mathcal J}\Lambda^{(j)}_0(p	^*)
\end{equation}
$\mathbb P$-a.s., where $\phi$ denotes a generic random characteristics which in distribution equals $\phi^{(j)}$, $j\in\mathcal J$. In order to  deduce (\ref{e.Prop_8}), note that Lemma~6 of \cite{Kno12} yields that
\begin{align}\label{e.L_6}
\lim_{\eta\downarrow0}\mathbb E\left(\eta^{1+p^*}Z^I_\eta\right) &= \lim_{\eta\downarrow0}\mathbb E\left(\eta^{1+p^*}Z^\phi_\eta\right)\mathbb E\left(\sum_{j\in\mathcal J}v_j^{1+p^*}\right)
\\[0.5ex]
&= \mathbb E\left(\frac{\Lambda_0(p^*)}{\Phi'(p^*)}\int_{(0,1)}\sum_{k\in\N}|\Pi_{k}(t)|^{1+p^*}\int_\mathcal P\int_{(0,1]}s^{p^*}\phi\left(\rho,\pi\right)\dd s\mu(\dd\pi)\dd t\right)\mathbb E\left(\sum_{j\in\mathcal J}v_j^{1+p^*} \right)\notag
\end{align}
and analogously, by resorting to Lemma~9 of \cite{Kno12}, we also obtain that the map $\eta\mapsto\mathbb E(\eta^{1+p^*}Z^I_\eta)$ is continuous on $(0,1]$. Hence, in view of  Lemma~\ref{l.L_3} we obtain that the analogue of Lemma~11 of \cite{Kno12} holds in the setting with immigration. By means of Proposition~\ref{l.t.SLLN.fprc.1b.new.1b.0} we can thus follow the lines of the proof of Proposition~7 of \cite{Kno12} to infer that (\ref{e.Prop_8}) holds.

Let us point out  that the key estimate in the proof of Theorem~1 of \cite{Kno12} to deduce the desired almost sure convergence, as asserted in that theorem, from Proposition~7 in the same reference is based on the following estimate (cf. (27) in \cite{Kno12}):
\begin{equation}\label{e.27}
\mathbb E\left(\rho^{\iota r(1+p^*)}\left|Z^{\phi}_{\rho^{\iota r}}-Z^{\phi_{\iota,\rho^r}}_{\rho^{\iota r}}\right|\right)
\le aA_r\rho^{-(1+p^*)}\mathbb E\left(\sup_{u\in(0,1]}\left(u^{1+p^*}T^{(a)}_{u}\right)\right)\frac{\rho^{(\iota-1)(p^*-\tilde p)r}}{1-\rho^{p^*-\tilde p}},
\end{equation}
where
\[
A_r:=\sup_{\eta\in(0,\rho^{(\iota-1)r}]}\int_\mathcal P\eta^{1+\tilde p}\mathbb E\left(\phi\left(\eta,\pi\right)\right)\mu(\dd\pi)
\]
and  $\tilde p\in(\underline p,p^*)$ is given by (\ref{e.assumption_1}). Since (\ref{e.27})  is only concerned with the distributions of $\Pi$ and $\phi$ and since we assumed that the immigrating processes as well as the associated random characteristics are identically distributed, we infer that in the setting with immigration the analogue of the estimate in (27) of \cite{Kno12} is uniform in $j\in\mathcal J$. Therefore, by means of (\ref{e.fpwi.1b_L1}) the convergence in (\ref{e.Prop_8}) can be extended to the desired almost sure convergence in (\ref{e.a.s.convergence}) along the lines of Parts~I and II of the proof of Theorem~1 in \cite{Kno12}. In the light of (\ref{e.L_6}) the $\mathscr L^1$-convergence in (\ref{e.a.s.convergence}) follows analogously to Part~ III in the proof of Theorem~1 of \cite{Kno12}. 
\hfill$\square$

\section{Rate of decay of the largest block}\label{s.decay}
This section is devoted to the proof of Theorem~\ref{t.1.2b}. In view of the martingale arguments used to prove similar results for fragmentation processes (cf. Corollary~1.4 of \cite{Ber06}) and fragmentation processes with killing (see Theorem~4 of \cite{KK12}) we aim at using  the submartingale convergence theorem to prove Theorem~\ref{t.1.2b}.  As in the previous two sections we consider throughout this section, unless stated otherwise, a deterministic vector $(v_j)_{j\in\mathcal J}\in\mathcal S$. The extension to the random vector $(v_j)_{j\in\mathcal J}$ we are interested in then holds, since Theorem~\ref{t.1.2b} is concerned with  almost sure properties.

Throughout this section we consider the set-up established in  Section~\ref{intro_immigration}. Moreover, we assume that the $\Pi^{(j)}$, $j\in\mathcal J$, are identically distributed and as before we omit the index $j$ for objects that are determined by the distribution of $\Pi^{(j)}$. 
 Note that this implies in particular that there exists some $\Phi$ with $\Phi^{(j)}=\Phi$ for all $j\in\mathcal J$. Recall that for every $j\in\mathcal J$  the time at which $\Pi^{(j)}$ immigrates is denoted by $\mathfrak t_j$. Furthermore, we define
\[
\mathfrak J_t:=\{j\in\mathcal J:\mathfrak t_j\le t\}.
\]
to be the set of all those indices belonging to blocks that immigrated not later than time $t\in\rpn$. For every $p\in(\underline p,\infty)$ consider the stochastic process $M^I(p):=(M^I_t(p))_{t\in\rpn}$ defined by
\[
\sum_{j\in\mathfrak J_t}v_j^{1+p}\sum_{n\in\N}|\Pi^{(j)}_n(t-\mathfrak t_j)|^{1+p}e^{\Phi(p)(t-\mathfrak t_j)}
\]
for all $t\in\rpn$

The following lemma shows that the martingale property of 
\[
M^{(j)}(p):=\left(\sum_{n\in\N}|\Pi^{(j)}_n(t)|^{1+p}e^{\Phi^{(j)}(p)t}\right)_{t\in\rpn},
\] 
cf. (\ref{e.martingale_M}), for all $j\in\mathcal J$ yields the submartingale property of $M^I(p)$.

\begin{lemma}\label{l.submartingale_M_I}
Assume that the $\Pi^{(j)}$, $j\in\mathcal J$, are identically distributed and that (\ref{e.fpwi.1b}) holds. Then $M^I(p)$ is an $\mathscr F^I$-submartingale for each $p\in(\underline p,\infty)$.
\end{lemma}

\begin{pf}
Note first that
\begin{align}\label{e.submartingale_decomposition.1}
\mathbb E\left(\left.M^I_{t+s}(p)\right|\mathscr F^I_t\right)
&=\mathbb E\left(\left.\sum_{j\in\mathfrak J_t}v_j^{1+p}\sum_{n\in\N}|\Pi^{(j)}_n(t+s-\mathfrak t_j)|^{1+p}e^{\Phi(p)(t+s-\mathfrak t_j)}\right|\mathscr F^I_t\right)
\\[0.5ex]
&\qquad +\mathbb E\left(\left.\sum_{j\in\mathfrak J_{t+s}\setminus\mathfrak J_t}v_j^{1+p}\sum_{n\in\N}|\Pi^{(j)}_n(t+s-\mathfrak t_j)|^{1+p}e^{\Phi(p)(t+s-\mathfrak t_j)}\right|\mathscr F^I_t\right).\notag
\end{align}
Let us now deal with the first summand on the right-hand side of (\ref{e.submartingale_decomposition.1}). To this end, observe that by means of the martingale property of $M^{(j)}(p)$ for every $j\in\mathcal J$ we have
\begin{align*}
&\mathbb E\left(\left.\sum_{j\in\mathfrak J_t}v_j^{1+p}\sum_{n\in\N}|\Pi^{(j)}_n(t+s-\mathfrak t_j)|^{1+p}e^{\Phi(p)(t+s-\mathfrak t_j)}\right|\mathscr F^I_t\right) 
\\[0.5ex]
&= \sum_{j\in\mathfrak J_t}v_j^{1+p}\left.\mathbb E\left(\left.\sum_{n\in\N}|\Pi^{(j)}_n(t+s-u)|^{1+p}e^{\Phi(p)(t+s-u)}\right|\mathscr F^{(j)}_{t-u}\right)\right|_{u=\mathfrak t_j}
\\[0.5ex]
&= \sum_{j\in\mathfrak J_t}v_j^{1+p}\sum_{n\in\N}|\Pi^{(j)}_n(t-\mathfrak t_j)|^{1+p}e^{\Phi(p)(t-\mathfrak t_j)}.
\end{align*}
Moreover, since the second summand on the right-hand side of (\ref{e.submartingale_decomposition.1}) is nonnegative, we thus infer that
\[
\mathbb E\left(\left.M^I_{t+s}(p)\right|\mathscr F_t\right)\ge\sum_{j\in\mathfrak J_t}v_j^{1+p}\sum_{n\in\N}|\Pi^{(j)}_n(t+s-\mathfrak t_j)|^{1+p}e^{\Phi(p)(t+s-\mathfrak t_j)}=M^I_t(p).
\]
Consequently, since $M^I(p)$ is an integrable and $\mathscr F^I$-adapted process, it is an $\mathscr F^I$-submartingale.
\end{pf}

In the proof of Theorem~\ref{t.1.2b} we shall use the following corollary of the previous lemma.
\begin{corollary}\label{c.pos_lim_submartingale}
Assume that the $\Pi^{(j)}$, $j\in\mathcal J$, are identically distributed and that (\ref{e.fpwi.1b}) holds. In addition, let $p\in(\underline p,\bar p)$. Then the limit $M^I_\infty(p):=\lim_{t\to\infty}M^I_t(p)$ exists and is positive $\mathbb P$-almost surely.
\end{corollary}

\begin{pf}
In view of Lemma~\ref{l.submartingale_M_I} the $\mathbb P$-a.s. existence of $M^I_\infty(p)$ follows from the submartingale convergence theorem and
\[
\sup_{t\in\rpn}\mathbb E\left(M^I_t\right) \le \mathbb E\left(\sum_{j\in\mathfrak J_t}v_j^{1+p}\left.\sup_{t\in\rpn}\mathbb E\left(\sum_{n\in\N}|\Pi^{(j)}_n(t-u)|^{1+p}e^{\Phi(p)(t-u)}\right)\right|_{u=\mathfrak t_j}\right)
\le \sum_{j\in\mathcal J}v_j^{1+p}
<\infty,
\]
where we used (\ref{e.fpwi.1b}). 
The $\mathbb P$-a.s positivity of $M^I_\infty(p)$ follows from 
\[
M^I_\infty(p)\ge v_{j_s}^{1+p}\lim_{t\to\infty}\sum_{n\in\N}|\Pi^{(j_s)}_n(t-\mathfrak t_{j_s})|^{1+p}\exp\left(\Phi^{(j_s)}(p)(t-\mathfrak t_{j_s})\right)>0
\]
for all $s\in\rpn$, cf. Theorem~2 of \cite{Ber03}.
\end{pf}

{\bf Proof of Theorem~\ref{t.1.2b}}
Note first that
 we have
\[
\left(\lambda^I_1(t)\right)^{1+\bar p}e^{\Phi(\bar p)(t-\mathfrak t_{j_t})} \le M^I_t(\bar p).
\]
holds for all $t\in\rpn$, where $j_t$ was defined in (\ref{e.j_t}).
Consequently, by taking the logarithm on both sides above, dividing  by $t$ and subsequently taking the limsup as $t\to\infty$, we deduce by resorting to Corollary~\ref{c.pos_lim_submartingale} that
\begin{equation}\label{e.t.1.2b.up}
\limsup_{t\to\infty}\frac{\ln\left(\lambda^I_1(t)\right)}{t}\le\left(\limsup_{t\to\infty}\frac{\mathfrak t_{j_t}}{t}-1\right)\frac{\Phi(\bar p)}{1+\bar p}=\left(\limsup_{t\to\infty}\frac{\mathfrak t_{j_t}}{t}-1\right)\Phi'(\bar p)
\end{equation}
holds $\mathbb P$-almost surely.

In order to show the converse inequality let $p\in(\underline p,\bar p)$ as well as $\epsilon\in(0,p-\underline p)$  and observe that
\begin{align*}
M^I_t(p) 
&\le \left(\lambda^I_1(t)\right)^\epsilon \sum_{j\in\mathfrak J_t}e^{(\Phi(p)-\Phi(p-\epsilon))(t-\mathfrak t_j)}v_j^{1+p-\epsilon}\sum_{n\in\N}\lambda^{(j)}_n(t-\mathfrak t_j)^{1+p-\epsilon}e^{\Phi(p-\epsilon)(t-\mathfrak t_j)}
\\[0.5ex]
&\le \left(\lambda^I_1(t)\right)^\epsilon e^{(\Phi(p)-\Phi(p-\epsilon))t}\sum_{j\in\mathfrak J_t}v_j^{1+p-\epsilon}\sum_{n\in\N}\lambda^{(j)}_n(t-\mathfrak t_j)^{1+p-\epsilon}e^{\Phi(p-\epsilon)(t-\mathfrak t_j)}
\\[0.5ex]
&\le \left(\lambda^I_1(t)\right)^\epsilon e^{(\Phi(p)-\Phi(p-\epsilon))t} M^I_t(p-\epsilon). 
\end{align*}
Similarly to above, this time applying Corollary~\ref{c.pos_lim_submartingale} to both submartingales $M^I(p)$ and $M^I(p-\epsilon)$, we obtain that
\[
\liminf_{t\to\infty}\frac{\ln\left(\lambda^I_1(t)\right)}{t}\ge-\frac{\Phi(p)-\Phi(p-\epsilon)}{\epsilon}
\]
$\mathbb P$-almost surely. Hence, letting $\epsilon\downarrow0$  results in
\[
\liminf_{t\to\infty}\frac{\ln\left(\lambda^I_1(t)\right)}{t}\ge-\lim_{\epsilon\downarrow0}\frac{\Phi(p)-\Phi(p-\epsilon)}{\epsilon}=-\Phi'(p)
\]
$\mathbb P$-almost surely.
By means of the continuity of $\phi'$ we thus infer by letting $p\uparrow\bar p$ that
\[
\liminf_{t\to\infty}\frac{\ln\left(\lambda^I_1(t)\right)}{t}\ge-\Phi'(\bar p)
\]
holds $\mathbb P$-almost surely. In view of (\ref{e.t.1.2b.up}) this completes the proof.
\hfill$\square$

\section*{Acknowledgement}
I would like to thank 
Andreas E. Kyprianou for the helpful discussions regarding the topic of this paper.

\end{document}

%% file: frag_stopping_line_50_eta.pstex_t
\begin{picture}(0,0)%
\includegraphics{frag_stopping_line_50_eta.pstex}%
\end{picture}%
\setlength{\unitlength}{2072sp}%
\begingroup\makeatletter\ifx\SetFigFont\undefined%
\gdef\SetFigFont#1#2#3#4#5{%
  \reset@font\fontsize{#1}{#2pt}%
  \fontfamily{#3}\fontseries{#4}\fontshape{#5}%
  \selectfont}%
\fi\endgroup%
\begin{picture}(7149,6781)(121,-6515)
\put(5581,-196){\makebox(0,0)[lb]{\smash{{\SetFigFont{6}{7.2}{\rmdefault}{\mddefault}{\updefault}{\color[rgb]{0,0,0}$\vdots$}%
}}}}
\put(6706,119){\makebox(0,0)[lb]{\smash{{\SetFigFont{6}{7.2}{\rmdefault}{\mddefault}{\updefault}{\color[rgb]{0,0,0}$\vdots$}%
}}}}
\put(6886,-6271){\makebox(0,0)[lb]{\smash{{\SetFigFont{6}{7.2}{\rmdefault}{\mddefault}{\updefault}{\color[rgb]{0,0,0}time}%
}}}}
\put(4456,-781){\makebox(0,0)[lb]{\smash{{\SetFigFont{6}{7.2}{\rmdefault}{\mddefault}{\updefault}{\color[rgb]{0,0,0}$\vdots$}%
}}}}
\put(3421,-331){\makebox(0,0)[lb]{\smash{{\SetFigFont{6}{7.2}{\rmdefault}{\mddefault}{\updefault}{\color[rgb]{0,0,0}$\vdots$}%
}}}}
\put(6841,-1591){\makebox(0,0)[lb]{\smash{{\SetFigFont{6}{7.2}{\rmdefault}{\mddefault}{\updefault}{\color[rgb]{0,0,0}$\ldots$}%
}}}}
\put(6841,-1276){\makebox(0,0)[lb]{\smash{{\SetFigFont{6}{7.2}{\rmdefault}{\mddefault}{\updefault}{\color[rgb]{0,0,0}$\ldots$}%
}}}}
\put(6841,-736){\makebox(0,0)[lb]{\smash{{\SetFigFont{6}{7.2}{\rmdefault}{\mddefault}{\updefault}{\color[rgb]{0,0,0}$\ldots$}%
}}}}
\put(6841,-61){\makebox(0,0)[lb]{\smash{{\SetFigFont{6}{7.2}{\rmdefault}{\mddefault}{\updefault}{\color[rgb]{0,0,0}$\ldots$}%
}}}}
\put(5716,-736){\makebox(0,0)[lb]{\smash{{\SetFigFont{6}{7.2}{\rmdefault}{\mddefault}{\updefault}{\color[rgb]{0,0,0}$\ldots$}%
}}}}
\put(5716,-376){\makebox(0,0)[lb]{\smash{{\SetFigFont{6}{7.2}{\rmdefault}{\mddefault}{\updefault}{\color[rgb]{0,0,0}$\ldots$}%
}}}}
\put(3556,-1501){\makebox(0,0)[lb]{\smash{{\SetFigFont{6}{7.2}{\rmdefault}{\mddefault}{\updefault}{\color[rgb]{0,0,0}$\ldots$}%
}}}}
\put(3556,-1186){\makebox(0,0)[lb]{\smash{{\SetFigFont{6}{7.2}{\rmdefault}{\mddefault}{\updefault}{\color[rgb]{0,0,0}$\ldots$}%
}}}}
\put(3556,-511){\makebox(0,0)[lb]{\smash{{\SetFigFont{6}{7.2}{\rmdefault}{\mddefault}{\updefault}{\color[rgb]{0,0,0}$\ldots$}%
}}}}
\put(1891,-5371){\makebox(0,0)[lb]{\smash{{\SetFigFont{6}{7.2}{\rmdefault}{\mddefault}{\updefault}{\color[rgb]{0,0,0}$\ldots$}%
}}}}
\put(1891,-4516){\makebox(0,0)[lb]{\smash{{\SetFigFont{6}{7.2}{\rmdefault}{\mddefault}{\updefault}{\color[rgb]{0,0,0}$\ldots$}%
}}}}
\put(1891,-3751){\makebox(0,0)[lb]{\smash{{\SetFigFont{6}{7.2}{\rmdefault}{\mddefault}{\updefault}{\color[rgb]{0,0,0}$\ldots$}%
}}}}
\put(2791,-1771){\makebox(0,0)[lb]{\smash{{\SetFigFont{6}{7.2}{\rmdefault}{\mddefault}{\updefault}{\color[rgb]{0,0,0}$\ldots$}%
}}}}
\put(2791,-2581){\makebox(0,0)[lb]{\smash{{\SetFigFont{6}{7.2}{\rmdefault}{\mddefault}{\updefault}{\color[rgb]{0,0,0}$\ldots$}%
}}}}
\put(2791,-2986){\makebox(0,0)[lb]{\smash{{\SetFigFont{6}{7.2}{\rmdefault}{\mddefault}{\updefault}{\color[rgb]{0,0,0}$\ldots$}%
}}}}
\put(2791,-3571){\makebox(0,0)[lb]{\smash{{\SetFigFont{6}{7.2}{\rmdefault}{\mddefault}{\updefault}{\color[rgb]{0,0,0}$\ldots$}%
}}}}
\put(3556,-1906){\makebox(0,0)[lb]{\smash{{\SetFigFont{6}{7.2}{\rmdefault}{\mddefault}{\updefault}{\color[rgb]{0,0,0}$\ldots$}%
}}}}
\put(4591,-1546){\makebox(0,0)[lb]{\smash{{\SetFigFont{6}{7.2}{\rmdefault}{\mddefault}{\updefault}{\color[rgb]{0,0,0}$\ldots$}%
}}}}
\put(4591,-1141){\makebox(0,0)[lb]{\smash{{\SetFigFont{6}{7.2}{\rmdefault}{\mddefault}{\updefault}{\color[rgb]{0,0,0}$\ldots$}%
}}}}
\put(4591,-961){\makebox(0,0)[lb]{\smash{{\SetFigFont{6}{7.2}{\rmdefault}{\mddefault}{\updefault}{\color[rgb]{0,0,0}$\ldots$}%
}}}}
\put(1756,-2581){\makebox(0,0)[lb]{\smash{{\SetFigFont{6}{7.2}{\rmdefault}{\mddefault}{\updefault}{\color[rgb]{0,0,0}$\vdots$}%
}}}}
\put(2656,-1546){\makebox(0,0)[lb]{\smash{{\SetFigFont{6}{7.2}{\rmdefault}{\mddefault}{\updefault}{\color[rgb]{0,0,0}$\vdots$}%
}}}}
\put(5716,-2446){\makebox(0,0)[lb]{\smash{{\SetFigFont{6}{7.2}{\rmdefault}{\mddefault}{\updefault}{\color[rgb]{0,0,0}$\ldots$}%
}}}}
\put(5716,-2761){\makebox(0,0)[lb]{\smash{{\SetFigFont{6}{7.2}{\rmdefault}{\mddefault}{\updefault}{\color[rgb]{0,0,0}$\ldots$}%
}}}}
\put(4591,-2356){\makebox(0,0)[lb]{\smash{{\SetFigFont{6}{7.2}{\rmdefault}{\mddefault}{\updefault}{\color[rgb]{0,0,0}$\ldots$}%
}}}}
\put(5716,-3346){\makebox(0,0)[lb]{\smash{{\SetFigFont{6}{7.2}{\rmdefault}{\mddefault}{\updefault}{\color[rgb]{0,0,0}$\ldots$}%
}}}}
\put(721,-6451){\makebox(0,0)[lb]{\smash{{\SetFigFont{6}{7.2}{\rmdefault}{\mddefault}{\updefault}{\color[rgb]{0,0,0}$0$}%
}}}}
\put(136,-3256){\makebox(0,0)[lb]{\smash{{\SetFigFont{6}{7.2}{\rmdefault}{\mddefault}{\updefault}{\color[rgb]{0,0,0}$-\ln(\eta)$}%
}}}}
\put(946,-196){\makebox(0,0)[lb]{\smash{{\SetFigFont{6}{7.2}{\rmdefault}{\mddefault}{\updefault}{\color[rgb]{0,0,0}$-\ln(\text{size})$}%
}}}}
\end{picture}%

%% file: frag_stopping_line_50b_eta.pstex_t
\begin{picture}(0,0)%
\includegraphics{frag_stopping_line_50b_eta.pstex}%
\end{picture}%
\setlength{\unitlength}{2072sp}%
\begingroup\makeatletter\ifx\SetFigFont\undefined%
\gdef\SetFigFont#1#2#3#4#5{%
  \reset@font\fontsize{#1}{#2pt}%
  \fontfamily{#3}\fontseries{#4}\fontshape{#5}%
  \selectfont}%
\fi\endgroup%
\begin{picture}(7149,6490)(121,-6515)
\put(5581,-196){\makebox(0,0)[lb]{\smash{{\SetFigFont{6}{7.2}{\rmdefault}{\mddefault}{\updefault}{\color[rgb]{0,0,0}$\vdots$}%
}}}}
\put(6886,-6271){\makebox(0,0)[lb]{\smash{{\SetFigFont{6}{7.2}{\rmdefault}{\mddefault}{\updefault}{\color[rgb]{0,0,0}time}%
}}}}
\put(4456,-781){\makebox(0,0)[lb]{\smash{{\SetFigFont{6}{7.2}{\rmdefault}{\mddefault}{\updefault}{\color[rgb]{0,0,0}$\vdots$}%
}}}}
\put(1891,-5371){\makebox(0,0)[lb]{\smash{{\SetFigFont{6}{7.2}{\rmdefault}{\mddefault}{\updefault}{\color[rgb]{0,0,0}$\ldots$}%
}}}}
\put(1891,-4516){\makebox(0,0)[lb]{\smash{{\SetFigFont{6}{7.2}{\rmdefault}{\mddefault}{\updefault}{\color[rgb]{0,0,0}$\ldots$}%
}}}}
\put(1891,-3751){\makebox(0,0)[lb]{\smash{{\SetFigFont{6}{7.2}{\rmdefault}{\mddefault}{\updefault}{\color[rgb]{0,0,0}$\ldots$}%
}}}}
\put(2791,-3571){\makebox(0,0)[lb]{\smash{{\SetFigFont{6}{7.2}{\rmdefault}{\mddefault}{\updefault}{\color[rgb]{0,0,0}$\ldots$}%
}}}}
\put(1756,-2581){\makebox(0,0)[lb]{\smash{{\SetFigFont{6}{7.2}{\rmdefault}{\mddefault}{\updefault}{\color[rgb]{0,0,0}$\vdots$}%
}}}}
\put(2656,-1546){\makebox(0,0)[lb]{\smash{{\SetFigFont{6}{7.2}{\rmdefault}{\mddefault}{\updefault}{\color[rgb]{0,0,0}$\vdots$}%
}}}}
\put(5716,-3346){\makebox(0,0)[lb]{\smash{{\SetFigFont{6}{7.2}{\rmdefault}{\mddefault}{\updefault}{\color[rgb]{0,0,0}$\ldots$}%
}}}}
\put(946,-196){\makebox(0,0)[lb]{\smash{{\SetFigFont{6}{7.2}{\rmdefault}{\mddefault}{\updefault}{\color[rgb]{0,0,0}$-\ln(\text{size})$}%
}}}}
\put(721,-6451){\makebox(0,0)[lb]{\smash{{\SetFigFont{6}{7.2}{\rmdefault}{\mddefault}{\updefault}{\color[rgb]{0,0,0}$0$}%
}}}}
\put(136,-3256){\makebox(0,0)[lb]{\smash{{\SetFigFont{6}{7.2}{\rmdefault}{\mddefault}{\updefault}{\color[rgb]{0,0,0}$-\ln(\eta)$}%
}}}}
\end{picture}%